\documentclass[12pt]{article}
\usepackage{amsmath}
\usepackage{graphicx}
\usepackage{enumerate}
\usepackage{natbib}
\usepackage{url} 

\RequirePackage{amsthm,amsmath,amsfonts,amssymb}

\usepackage{hyperref}

\newcommand{\blind}{1}

\usepackage[ruled]{algorithm2e}

\usepackage{bm}

\makeatletter
\newcommand{\rmnum}[1]{\romannumeral #1}
\newcommand{\Rmnum}[1]{\expandafter\@slowromancap\romannumeral #1@}
\makeatother

\usepackage{tabu}                     
\usepackage{multirow}                 
\usepackage{multicol}                 
\usepackage{multirow}                
\usepackage{float}                    
\usepackage{makecell}                 
\usepackage{booktabs} 
\usepackage{subcaption}
\usepackage{xcolor}
\usepackage[figuresleft]{rotating}
\usepackage{pdflscape}

\theoremstyle{definition}
\newtheorem{Co}{Corollary}
\newtheorem{D}{Definition}

\newtheorem{ass}{Assumption}

\newtheorem*{remark}{Remark}

\newtheorem{manualtheoreminner}{Assumption}
\newenvironment{manualtheorem}[1]{%
  \renewcommand\themanualtheoreminner{#1}%
  \manualtheoreminner
}{\endmanualtheoreminner}

\def \tr{\text{tr}}

\def \op{\text{op}}

\makeatletter
\newcommand*{\addFileDependency}[1]{
  \typeout{(#1)}
  \@addtofilelist{#1}
  \IfFileExists{#1}{}{\typeout{No file #1.}}
}
\makeatother

\theoremstyle{plain}

\newtheorem{theorem}{Theorem}[section]

\theoremstyle{remark}

\begin{document}

\def\spacingset#1{\renewcommand{\baselinestretch}%
{#1}\small\normalsize} \spacingset{1}


\if1\blind
{
  \title{\bf Universal Bootstrap for Spectral Statistics: Beyond Gaussian Approximation}
  \author{Guoyu Zhang \thanks{
    guoyz@stu.pku.edu.cn.}\hspace{.2cm}\\
    Department of Probability and Statistics School of Mathematical Sciences,\\
    Center for Statistical Science,
Peking University\\
    Dandan Jiang \thanks{
    jiangdd@xjtu.edu.cn.}\hspace{.2cm}\\
    School of Mathematics and Statistics,
Xi'an Jiaotong University \\
    Fang Yao \thanks{
     Corresponding author, fyao@math.pku.edu.cn. This research is partially supported by the National Key Research and Development Program of China (No. 2022YFA1003801), the National Natural Science Foundation of China (No. 12292981, 12288101), the New Cornerstone Science Foundation through Xplorer Prize, the Fundamental Research Funds for the Central Universities (LMEQF), and the LMAM.
     }\hspace{.2cm}\\
    School of Mathematical Sciences,\\
    Center for Statistical Science,
Peking University}
\date{}
  \maketitle
} \fi

\if0\blind
{
  \bigskip
  \bigskip
  \bigskip
  \begin{center}
    {\LARGE\bf Universal Bootstrap for Spectral Statistics: Beyond Gaussian Approximation}
\end{center}
  \medskip
} \fi

\bigskip
\begin{abstract}
Spectral analysis plays a crucial role in high-dimensional statistics, where determining the asymptotic distribution of various spectral statistics remains a challenging task. Due
to the difficulties of deriving the analytic form, recent advances have explored data-driven bootstrap methods for this purpose. However, widely used Gaussian approximation-based bootstrap methods, such as the empirical bootstrap and multiplier bootstrap, have been shown to be inconsistent in approximating the distributions of spectral statistics in high-dimensional settings. To address this issue, we propose a universal bootstrap procedure based on the concept of universality from random matrix theory. Our method consistently approximates a broad class of spectral statistics across both high- and ultra-high-dimensional regimes, accommodating scenarios where the dimension-to-sample-size ratio $p/n$ converges to a nonzero constant or diverges to infinity without requiring structural assumptions on the population covariance matrix, such as eigenvalue decay or low effective rank. 
We showcase this universal bootstrap method for high-dimensional covariance inference. Extensive simulations and a real-world data study support our findings, highlighting the favorable finite sample performance of the proposed universal bootstrap procedure.
\end{abstract}

\noindent%
\vfill

\newpage
\spacingset{1.9} 
\section{Introduction}

\label{intro}

Spectral analysis of large covariance matrices is crucial in multivariate statistical estimation and hypothesis testing.
Consider $n$ independent, not necessarily identically distributed random vectors $\bm{X}_1,\cdots,\bm{X}_n \in \mathbb{R}^p$ with zero mean and a common covariance matrix $\bm{\Sigma}$. A prominent example we consider is the spectral statistic
\begin{align}\label{statistic}
    T = \left\|\bm{\hat{\Sigma}}- \bm{\Sigma}_0\right\|_{\text{op}},
\end{align}
where $\bm{\hat{\Sigma}}=\sum_{i=1}^n \bm{X}_i^{}\bm{X}_i^T/n$ is the sample covariance matrix, $\bm{\Sigma}_0$ is the pre-specified hypothesized covariance matrix, and the operator norm of a matrix $\bm{A}$ is defined as $\|\bm{A}\|_{\text{op}}=\sup_{\|\bm{u}\|=1}\| \bm{A}\bm{u}\|$. 
The statistic $T$ is fundamental in high-dimensional statistics, including covariance estimation \citep{fan2011high}, principal component analysis \citep{lopes2023bootstrapping}, and common factor determination \citep{yu2024testing}. Consequently, many studies have established non-asymptotic upper bounds on the tail behavior of $T$ \citep[e.g.,][]{adamczak2011sharp, bunea2015sample, koltchinskii2017concentration}. However, these bounds tend to be conservative when constructing confidence intervals for $T$. Few studies have examined the asymptotic distribution and bootstrap approximation of $T$ for general $\bm{\Sigma}$, $\bm{\Sigma}_0$ \citep{han2018gaussian, lopes2022improved, lopes2023bootstrapping, giessing2023gaussian}. Nonetheless, these works assume that either the dimension $p$ grows at a slower rate than $n$, or that the eigenvalues of $\bm{\Sigma}$ and $\bm{\Sigma}_0$ decay fast enough, making its effective rank smaller than $n$.
These assumptions enable Gaussian approximation but effectively lower the problem’s intrinsic dimension, leaving true high-dimensional cases unaddressed.

This article studies the behavior of $T$ and other general spectral statistics under both high- and ultra-high-dimensional regimes, allowing the dimension-to-sample size ratio $\phi = p/n$ to converge to a nonzero constant or diverge to infinity. Unlike previous studies, we impose no eigen-decay assumptions on $\bm{\Sigma}$ or $\bm{\Sigma}_0$.
Under these conditions, standard high-dimensional Gaussian approximations become invalid. To address this gap, we establish a new universality principle for spectral statistics, which leads to a theoretically justified and computationally efficient bootstrap method for approximating their distributions.

\subsection{A universality approach: overcoming Gaussianity}
The bootstrap method is a generic resampling technique to approximate the distribution of a statistic.  Originally designed for fixed-dimensional problems, the bootstrap method has recently been adapted for high-dimensional settings through substantial work. \citet{chernozhukov2013gaussian,chernozhukov2017central,deng2020beyond,lopes2020bootstrapping,lopes2022central,chernozhukov2023nearly} established Gaussian approximation rates for empirical and multiplier bootstrap for the supremum norm of high-dimensional mean vectors. Moreover, \citet{han2018gaussian,lopes2019bootstrapping,yao2021rates,lopes2022improved,lopes2023bootstrapping} explored bootstrap methods for spectral statistics of covariance matrices, typically assuming a low intrinsic dimension or effective rank relative to the sample size $n$. 

Despite these advances, \citet{el2019non,dette2024nonparametric} and \citet{yu2024testing} demonstrated that the Gaussian approximation-based bootstrap for spectral statistics fails when the intrinsic dimension is comparable to the sample size. 
These bootstrap methods approximate the distribution of the sample covariance matrix $\bm{\hat{\Sigma}}=\sum_{i=1}^n \bm{X}_i^{}\bm{X}_i^T/n$ using a Gaussian matrix with the same covariance. Nonetheless, as the dimension $p$ grows, this approximation becomes less accurate and eventually inconsistent, as shown by \citet{el2019non}. This suggests that high-dimensional covariance matrix structures constrain the applicability of the central limit theorem.

To overcome this limitation, we leverage the universality property of high-dimensional covariance structures.
Universality suggests that the asymptotic distribution of $T$ only depends on the distribution of $\bm{X}_1,\cdots,\bm{X}_n$ through their first two moments. This allows us to construct the universal bootstrap statistic by substituting independent Gaussian samples $\bm{Y}_1, \dots, \bm{Y}_n \sim \mathcal{N}(\bm{0},\bm{\Sigma})$ in place of $\bm{X}_1, \dots, \bm{X}_n$ in  the definition of (\ref{statistic}),
\begin{align}\label{ubstatistic}
    T^{\text{ub}}=\left\| \bm{\hat{\Sigma}}^{\text{ub}} -\bm{\Sigma}_0\right\|_{\text{op}},
\end{align}
where $\bm{\hat{\Sigma}}^{\text{ub}}=\sum_{i=1}^n \bm{Y}_i^{}\bm{Y}_i^T/n$.
The key insight is that although the bootstrapped matrix $\bm{\hat{\Sigma}}^{\text{ub}}$ does not approximate a Gaussian matrix, it effectively preserves the structure of $\bm{\hat{\Sigma}}$,
which empirical and multiplier bootstraps fail to maintain.
Despite its reliance on intricate random matrix theory, this method remains computationally straightforward to implement.
Further extensions for general spectral statistics beyond $T$ are also developed.

\subsection{A motivating example: covariance inference}
Our universal bootstrap method applies directly to spectral statistics in high-dimensional covariance tests. Testing large covariance has received considerable attention due to its critical role in high-dimensional statistics. Typically, hypothesis tests for high-dimensional covariance matrices $H_0: \bm{\Sigma}=\bm{\Sigma}_0$ rely on the norm between $\bm{\hat{\Sigma}}$ and $\bm{\Sigma}_0$, commonly the Frobenius norm \citep[e.g.,][]{chen2010tests, cai2013optimal} or the supremum norm  \citep[e.g.,][]{jiang2004asymptotic,cai2011limiting}, paralleling their high-dimensional mean tests \citet{chen2010two,tony2014two}. 
In contrast to the Frobenius norm and the supremum norm, which treat the covariance matrices as vectors and neglect their spectral structures, the operator norm captures the spectral structure and has recently gained significant attention. 
Various approaches have been developed for analyzing $T$ and its variants in specific scenarios for covariance tests.
When $\bm{\Sigma}_0=\bm{I}_p$, the distribution of $T$ can be derived by the extreme eigenvalues of the sample covariance matrix $\bm{\hat{\Sigma}}$. \citet{el2007tracy,bao2015universality,lee2016tracy,knowles2017anisotropic} established the limiting Tracy–Widom law for top eigenvalue of $\bm{\hat{\Sigma}}$ when dimension $p$ and sample size $n$ grow proportionally.
In cases where $\bm{\Sigma}_0$ is invertible, \citet{bao2015universality} proposed transforming the data $\bm{X}_1,\cdots,\bm{X}_n $ into $\bm{\Sigma}_0^{-1/2}\bm{X}_1,\cdots,\bm{\Sigma}_0^{-1/2}\bm{X}_n$ and considered the Roy’s largest root statistic
\begin{align}\label{roy}
    T^{\text{Roy}} = \left\|\bm{\Sigma}_0^{-\frac{1}{2}}\bm{\hat{\Sigma}}\bm{\Sigma}_0^{-\frac{1}{2}}- \bm{I}_p\right\|_{\text{op}}.
\end{align}
The limiting Tracy–Widom distribution also holds for $T^{\text{Roy}}$. However, for general $\bm{\Sigma}_0$, distributional results for $T$ remain unknown, even when $p/n$ is bounded. The universal bootstrap method addresses this gap by leveraging universality in random matrix theory. 
More applications in high-dimensional statistics are presented in Section \ref{apps}.

\subsection{Our contributions}
We summarize our contribution as three-fold. 
First, we propose a universal bootstrap procedure for $T$ based on universality results. We show that this procedure effectively approximates distributions of spectral statistics. As a byproduct, we establish the Tracy–Widom law for the largest eigenvalue of sample covariance matrices with general entries in the ultra-high-dimensional regime, addressing a long-standing gap since \citet{karoui2003largest}.
Additionally, the asymptotic distribution of $T$ is shown to be independent of the third and fourth moments of entries. In contrast, the Frobenius norm and supremum norm lack this universality.
Consequently, statistics based on Frobenius and supremum norms require fourth-moment estimation, which is generally more complex and computationally intensive. This underscores a key advantage of using operator norm-based spectral statistics.

Second,  we develop a unified bootstrap method for general spectral statistics and apply it to the covariance inference. 
We conduct power analysis for $T$ and $T^{\text{Roy}}$ within the generalized spiked model framework from \citet{jiang2021generalized}. We show that $T^{\text{Roy}}$ performs better in worst-case scenarios while $T$ excels in average performance. To further enhance the power, we propose a new combined statistic, $T^{\text{Com}}$, supported by a generalized universality result that serves as a theoretical guarantee for the universal bootstrap applicable to general spectral statistics.
Extensive simulations validate these findings, highlighting the superior performance of our combined statistics across a wide range of scenarios.

Third, from a technical perspective, we establish the anisotropic local law and universality for matrices of the form $\hat{\bm{\Sigma}}+\bm{R}$ in both high- and ultra-high-dimensional settings. When $\bm{R}=-\bm{\Sigma}_0$, this result characterizes the limiting properties of $T$. The anisotropic local law for $\bm{R}=\bm{0}$ was previously established by \citet{knowles2017anisotropic} for bounded $p/n$ and by \citet{ding2023global, ding2024two} for diverging $p/n$. However, the presence of $\bm{R}$ changes the problem structure, and the ultra-high-dimensional setting makes existing bounds suboptimal. We address this gap by introducing a new weighted anisotropic local law, refining the results of \citet{ding2023global} and \citet{ding2024two}. This advancement strengthens the framework for establishing universality.

\subsection{Notations and paper organization}
Throughout the paper, we reserve boldfaced symbols for vectors and matrices.
For a complex number $z\in\mathbb{C}$, we use $\Re z$ and $\Im z$ for its real part and imaginary part, respectively.
For a vector $\bm{u}\in \mathbb{R}^p$, we use $\|\bm{u}\|=\sqrt{\sum_{i=1}^p u_i^2}$ for its Euclidean norm. 
Denote the singular values of a matrix $\bm{A}$ by $\sigma_1(\bm{A})\geq\sigma_2(\bm{A})\geq \cdots \geq\sigma_r(\bm{A})$, where $r=\text{rank}(\bm{A})$. We use $\lambda_1(\bm{A}),\lambda_2(\bm{A}), \cdots \lambda_M(\bm{A})$ for the eigenvalues of a squared $\bm{A}$. 
For two sequences $\{a_n\}_{n=1}^\infty$, $\{b_n\}_{n=1}^\infty$, we use $a_n \lesssim b_n$ or $a_n=O(b_n)$ to show there exists a constant $C$ not depending  on $n$ such that $|a_n| \leq C b_n$ for all $n\in\mathbb{N}$. Denote $a_n=o(b_n)$ if $a_n/b_n\to 0$. 

The paper is organized as follows. Section \ref{mainre} provides an overview of the proposed method. Section \ref{univ} presents the main theoretical results on the consistency of the universal bootstrap procedure.
Section \ref{apps} discusses the applications of our universal bootstrap in the covariance test and construction of simultaneous confidence intervals.
Section \ref{numericalres} provides simulation results and a real data example, demonstrating the numerical performance of our methods. Detailed technical proofs, as well as the codes and datasets for reproducibility, are provided in the supplementary material. 

\section{Universal bootstrap}\label{mainre}
In this section, we present the universal bootstrap procedure and its generalized version.
Consider $n$ independent random vectors $\bm{X}_1,\cdots,\bm{X}_n \in \mathbb{R}^p$ with zero mean and covariance matrix $\bm{\Sigma}=(\sigma_{ij})_{p\times p}$, which are {\em not necessarily} identically distributed. 
As we focus on the covariance structure, we allow the third and fourth moments of $\bm{X}_i$ to differ across $i=1,\cdots,n$.
To simplify notation, we arrange the data into an $n \times p$ matrix $\bm{X}=(\bm{X}_1,\cdots,\bm{X}_n)^T \in \mathbb{R}^{n \times p}$, where $\bm{A}^T$ denotes the transpose of $\bm{A}$, and the sample covariance matrix is expressed by  $\bm{\hat{\Sigma}} = \bm{X}^T \bm{X}/n$.

Recalling the definition of statistic $T$ in (\ref{statistic}), we aim to characterize the asymptotic distribution of $T$. Consider the Gaussian data matrix $\bm{Y}=(\bm{Y}_1,\cdots,\bm{Y}_n)^T\in \mathbb{R}^{n\times p}$, with $\bm{Y}_1,\cdots,\bm{Y}_n \sim \mathcal{N}(\bm{0},\bm{\Sigma})$. Its sample covariance is defined accordingly as $\bm{\hat{\Sigma}}^{\text{ub}}=\bm{Y}^T\bm{Y}/n$.
In Section \ref{univ}, we demonstrate that the asymptotic distribution of $T$ can be uniformly approximated by $T^{\text{ub}}$, as defined in (\ref{ubstatistic}).
This motivates the universal bootstrap procedure. In practice, due to the analytical complexity of the distribution of $T^{\text{ub}}$, we generate $B$ independent samples, $\bm{Y}^1, \dots, \bm{Y}^B$, drawn from the same distribution as  $\bm{Y}$. We compute the sample covariance matrix $\bm{\hat{\Sigma}}^{\text{ub},b} = (\bm{Y}^b)^T(\bm{Y}^b)/n$ for each $b = 1, \dots, B$ and define
\begin{equation*}
    T^{\text{ub},b}=\left\| \bm{\hat{\Sigma}}^{\text{ub},b} -\bm{\Sigma}_0\right\|_{\text{op}}, \ b=1,\cdots,B.
\end{equation*}
The empirical distribution of $T^{\text{ub},b}$ for $b = 1, \dots, B$ serves as an approximation for the distribution of $T^{\text{ub}}$. Specifically, we use the empirical upper-$\alpha$ quantile, $\hat{q}^{\text{ub},B}_{\bm{\Sigma},\bm{\Sigma}_0}(\alpha)$, of $T^{\text{ub},b}$, $b=1, \dots, B$, as the bootstrap approximated upper-$\alpha$ quantile of $T$. In particular, when testing $H_0: \bm{\Sigma}=\bm{\Sigma}_0$, we will reject $H_0$ if $T>\hat{q}^{\text{ub},B}_{\bm{\Sigma}_0,\bm{\Sigma}_0}(\alpha)$, which does not contain unknown quantities.
The complete universal bootstrap algorithm is presented in Algorithm \ref{ubalgorithm}.

\begin{algorithm}
\caption{Universal bootstrap}\label{ubalgorithm}
  \SetKwInOut{Input}{input}\SetKwInOut{Output}{output}
\Input{Covariance matrix $\bm{\Sigma}$, $\bm{\Sigma}_0$, sample size $n$, test level $\alpha$, bootstrap number $B$}
  \Output{Bootstrapped $\alpha$ upper-quantile $\hat{q}^{\text{ub},B}_{\bm{\Sigma},\bm{\Sigma}_0}(\alpha)$}
  \For{$b\leftarrow 1$ \KwTo $B$}{
  $\bm{Z}^b \leftarrow    $ $n\times p$ matrix with i.i.d. standard normal entries \;
  $\bm{Y}^b \leftarrow  \bm{Z}^b \bm{\Sigma}^{1/2}  $  \;
  $\bm{\hat{\Sigma}}^{\text{ub},b} \leftarrow  (\bm{Y}^b)^T(\bm{Y}^b)/n  $  \;
  $T^{\text{ub},b} \leftarrow  \left\| \bm{\hat{\Sigma}}^{\text{ub},b} -\bm{\Sigma}_0\right\|_{\text{op}}  $  \;
  }
  $\hat{q}^{\text{ub},B}_{\bm{\Sigma},\bm{\Sigma}_0}(\alpha) \leftarrow $ empirical upper-$\alpha$ quantile of $T^{\text{ub},b}$, $b=1, \dots, B$.
\end{algorithm}

We now extend the universal bootstrap procedure for $T$ to general spectral statistics for the covariance test.
Similar to $T$ and $T^{\text{Roy}}$, such general spectral statistics can be constructed using the extreme singular values of $ \bm{\Sigma}_{2}^{-1/2}(\bm{\hat{\Sigma}}-\bm{\Sigma}_{1}) \bm{\Sigma}_{2}^{-1/2}$ for symmetric matrix $\bm{\Sigma}_{1}$ and positive-definite matrix $\bm{\Sigma}_{2}$. This inspires us to consider the statistic $T^{\text{ExS}}$ as a more general function of extreme singular values
\begin{align}\label{extrf}
    T^{\text{ExS}}=f\bigg(\bigg\{\left[\sigma_k\left( \bm{\Sigma}_{2,1}^{-\frac{1}{2}}(\bm{\hat{\Sigma}}-\bm{\Sigma}_{1,1}) \bm{\Sigma}_{2,1}^{-\frac{1}{2}}\right)\right]_{k=1}^{k_1},   \cdots,\left[\sigma_{k}\left( \bm{\Sigma}_{2,M}^{-\frac{1}{2}}(\bm{\hat{\Sigma}}-\bm{\Sigma}_{1,M}) \bm{\Sigma}_{2,M}^{-\frac{1}{2}}\right)\right]_{k=1}^{k_M}  \bigg\}\bigg),
\end{align}
where $f$ is a measurable function, $(\bm{\Sigma}_{1,m},\bm{\Sigma}_{2,m})$ are matrices, and $k_m$ are integers for $m=1,\cdots,M$. The statistic $T^{\text{ExS}}$ forms a general class of spectral statistics based on extreme singular values. For example, when $M=1$, $k_1=1$, $\bm{\Sigma}_{2,1}=\bm{I}_p$ and $\bm{\Sigma}_{1,1}=\bm{\Sigma}_0$, the statistic reduces to $T$. When $M=1$, $k_1=1$, $\bm{\Sigma}_{2,1}=\bm{\Sigma}_0$ and $\bm{\Sigma}_{1,1}=\bm{I}_p$, it reduces to $T^{\text{Roy}}$. This class also accommodates combinations of $T$ and $T^{\text{Roy}}$, regardless of their complex dependence. As discussed in Section \ref{apps}, $T^{\text{Roy}}$ is more effective in worst-case scenarios, while $T$ performs better on average.
To enhance testing power, we propose the following combined statistic
\begin{align}\label{comb}
    T^{\text{Com}}=\frac{T^2}{\text{tr}(\bm{\Sigma}_0)}+(T^{\text{Roy}})^2.
\end{align}
This statistic belongs to the general class of $T^{\text{ExS}}$.
Furthermore, while $T$ and $T^{\text{Roy}}$ rely only on the largest singular value (equivalent to the operator norm), $T^{\text{ExS}}$ incorporates the largest $k_m$ singular values, potentially improving test power \citet{ke2016detecting}.

The universal bootstrap extends naturally to the class of spectral statistics $T^{\text{ExS}}$.
The universal bootstrapped version, $T^{\text{ExS,ub}}$, is obtained by replacing $\bm{\hat{\Sigma}}$ with $\bm{\hat{\Sigma}}^{\text{ub}}$ in (\ref{extrf}). The empirical universal bootstrap threshold $\hat{q}^{\text{ExS,ub,B}}_{\bm{\Sigma},\left\{\bm{\Sigma}_{1,m},\bm{\Sigma}_{2,m}\right\}_{m=1}^{M}}(\alpha)$,  is defined as the upper $\alpha$-th quantile of $T^{\text{ExS,ub,1}},\cdots,T^{\text{ExS,ub,B}}$ of the i.i.d. samples $T^{\text{ExS,ub}}$. The algorithm for this generalized universal bootstrap is provided in Algorithm \ref{gubalgorithm}. In Section \ref{univ}, we establish the consistency of the universal bootstrap for generalized spectral statistics, demonstrating its applicability for constructing statistics using extreme eigenvalues in various ways.

\begin{algorithm}
\caption{Generalized universal bootstrap}\label{gubalgorithm}
  \SetKwInOut{Input}{input}\SetKwInOut{Output}{output}
\Input{Statistic function $f$, covariance matrix $\bm{\Sigma}$, $\left\{\bm{\Sigma}_{1,m},\bm{\Sigma}_{2,m}\right\}_{m=1}^{M}$ with number $k_1,\cdots,k_M$, sample size $n$, test level $\alpha$, bootstrap number $B$}
  \Output{Bootstrapped $\alpha$ upper-quantile $\hat{q}^{\text{ExS,ub,B}}_{\bm{\Sigma},\left\{\bm{\Sigma}_{1,m},\bm{\Sigma}_{2,m}\right\}_{m=1}^{M}}(\alpha)$}
  \For{$b\leftarrow 1$ \KwTo $B$}{
  $\bm{Z}^b \leftarrow    $ $n\times p$ matrix with i.i.d. standard normal entries \;
  $\bm{Y}^b \leftarrow  \bm{Z}^b \bm{\Sigma}^{1/2}  $  \;
  $\bm{\hat{\Sigma}}^{\text{ub},b} \leftarrow  (\bm{Y}^b)^T(\bm{Y}^b)/n  $  \;
  $T^{\text{ub},b} \leftarrow  f\bigg(\bigg\{\left[\sigma_k\left( \bm{\Sigma}_{2,1}^{-\frac{1}{2}}(\bm{\hat{\Sigma}}^{\text{ub},b}-\bm{\Sigma}_{1,1}) \bm{\Sigma}_{2,1}^{-\frac{1}{2}}\right)\right]_{k=1}^{k_1},   \cdots,\left[\sigma_{k}\left( \bm{\Sigma}_{2,M}^{-\frac{1}{2}}(\bm{\hat{\Sigma}}^{\text{ub},b}-\bm{\Sigma}_{1,M}) \bm{\Sigma}_{2,M}^{-\frac{1}{2}}\right)\right]_{k=1}^{k_M}  \bigg\}\bigg) $;
  }
  $\hat{q}^{\text{ExS,ub,B}}_{\bm{\Sigma},\left\{\bm{\Sigma}_{1,m},\bm{\Sigma}_{2,m}\right\}_{m=1}^{M}}(\alpha) \leftarrow $ empirical upper-$\alpha$ quantile of $T^{\text{ExS,ub,b}}$, $b=1, \dots, B$.
\end{algorithm}

\section{Theoretical guarantees}\label{univ}

\subsection{Universality properties}
\label{univer}
Our bootstrap consistency results rely on the universality from the random matrix theory. In this section, we 
review the background of random matrix theory and
establish the universality property.
To proceed, we first introduce some assumptions. We assume there exists some constants $C_1,C_2$ and $\gamma\geq 1$ such that 
\begin{align}\label{dimre}
    C_1n^{\gamma}\leq p \leq C_2 n^{\gamma}.
\end{align}
When $\gamma=1$, this is the typical proportional regime. When $\gamma>1$, the dimension-to-sample size ratio $p/n$ can diverge.
The primary matrix of interest is $\bm{\hat{\Sigma}}- \bm{\Sigma}_0$.
This inspires us to consider the matrix of the form $\bm{M}_n=\bm{\hat{\Sigma}}+\bm{R}$.
Here, $\bm{R}$ is a symmetric matrix that may have positive, negative, or zero eigenvalues. 
We impose some assumptions below.
\begin{ass}\label{ass1}
    Suppose that $\bm{X}_i=\bm{\Sigma}^{1/2}\bm{Z}_i$ and $\bm{Z}_i=(Z_{i1},\cdots,Z_{ip})^T$, and $Z_{ij}$, $ i=1,\cdots,n$, $j=1,\cdots,p $ are independent with $\mathbb{E}[Z_{ij}]=0$, $\mathbb{E}[Z_{ij}^2]=1$. There exists a positive sequence $C_k$ such that $\mathbb{E}[|Z_{ij}|^k]\leq C_k$ for $i=1,\cdots,n; j=1,\cdots,p$ and $k\in\mathbb{N}$.
\end{ass}

\begin{ass}\label{ass2}
    The matrix $\bm{\Sigma}$ and $\bm{R}$ are bounded in spectral norm, i.e. there exists some positive $C$ such that $\|\bm{\Sigma}\|_{\op},\ \|\bm{R}\|_{\op} \leq C$.
    Furthermore, there are constants $c_1, c_2 > 0$ such that the empirical spectral distribution $F$ of $\bm{\Sigma}$ satisfies $F^{\bm{\Sigma}}(c_1)\leq c_2$.
\end{ass}

\begin{ass}\label{ass3}
    The matrix $\bm{\Sigma}$ and $\bm{R}$ are commutative, i.e. $\bm{\Sigma}\bm{R}=\bm{R}\bm{\Sigma}$.
\end{ass}

Assumptions \ref{ass1} and \ref{ass2} are standard assumptions in the random matrix literature (\citet{,alex2014isotropic}). Notably, Assumption \ref{ass1} relies only on independence, without assuming identical distribution as in \citet{qiu2023asymptotic}. While all moments exist under Assumption \ref{ass1}, this condition could be relaxed (\citet{ding2018necessary}). We do not pursue this here. We also permit $\bm{\Sigma}$ to be singular, which is less restrictive than the invertibility assumed in \citet{ding2023global,ding2024two}.
Lastly, Assumption \ref{ass3}, imposed for technical requirements, necessitates that $\bm{\Sigma}$ and $\bm{R}$ share the same eigenvectors. The same assumption is required in the signal-plus-noise model as in \citet{zhou2023limiting,zhou2024analysis}.

We introduce a parameter $z$-dependent covariance $\bm{\Sigma}(z)=z\bm{\Sigma}(z\bm{I}_p-\phi^{-1/2}\bm{R})^{-1}$ and $\tilde{m}(z)$ as the fixed point of the following $z$-dependent deformed Mar\v{c}enko-Pastur law,
\begin{align}\label{zmp}
        \frac{1}{\tilde{m}(z)}=-z+\frac{1}{p}\sum_{i=1}^p \frac{\phi^{\frac{1}{2}} \lambda_i(\bm{\Sigma}(z))  }{1+\phi^{-\frac{1}{2}}\tilde{m}(z)\lambda_i(\bm{\Sigma}(z))}.
    \end{align}
This result mirrors the deformed Mar\v{c}enko-Pastur law in \citet{ding2023global}, with $\bm{\Sigma}(z)$ replacing the covariance matrix. The uniqueness and existence of the solution to (\ref{zmp}) are provided in the supplementary material. 
Define $E_{+}$ and $E_{-}$ to be the deterministic limit of $\lambda_1(\phi^{-1/2}\bm{M}_n)$ and $\lambda_p(\phi^{-1/2}\bm{M}_n)$, respectively, whose formal definition is provided in the supplementary material.
For the denominator in (\ref{zmp}), we impose the following technical assumption.
\begin{ass}\label{ass4}
    When $\gamma=1$ in (\ref{dimre}), we require that there exists $\tau>0$ such that
    \begin{align}\label{nonsp}
        | 1+\phi^{-\frac{1}{2}}\tilde{m}(E_{\pm}) \lambda_i(\bm{\Sigma}(E_{\pm}))| \geq \tau, \quad i=1,\cdots,p.
    \end{align}
\end{ass}
Informally, Assumption (\ref{nonsp}) ensures that the extreme eigenvalues of $\bm{\Sigma}$ do not spread near the endpoints $E_{\pm}$.
Similar assumptions have been made in the literature for the universality of $\hat{\bm{\Sigma}}$ \citep{bao2015universality,knowles2017anisotropic}, with the Stieltjes transform and $\lambda_i(\bm{\Sigma})$ replaced by $\tilde{m}$ and $\lambda_i(\bm{\Sigma}(E_{\pm}))$, respectively.  

We denote by $\mathbb{P}^{\text{Gau}}$ and $\mathbb{E}^{\text{Gau}}$ the probability and expectation under which $\{Z_{ij}\}_{1\leq i \leq n, 1\leq j\leq p}$ are independent standard Gaussian variables. With this, we present the following result.

\begin{theorem}[Universality of the largest eigenvalue]\label{gauappbd}
Under Assumptions \ref{ass1}, \ref{ass2}, \ref{ass3}, \ref{ass4}, there exists constant $C$ for large enough $n$ and $t\in\mathbb{R}$ such that for any small $\epsilon>0$, 
\begin{align}\label{gpbd}
    \mathbb{P}(n^{\frac{2}{3}} (&\lambda_1(\phi^{-1/2}\bm{M}_n)-E_+) \leq t-n^{-\epsilon} )-n^{-\frac{1}{6}+C\epsilon}\nonumber\\ 
    \leq	&\mathbb{P}^{\text{Gau}}(n^{\frac{2}{3}} (\lambda_1(\phi^{-1/2}\bm{M}_n)-E_+) \leq t )\\
    &\leq	\mathbb{P}(n^{\frac{2}{3}} (\lambda_1(\phi^{-1/2}\bm{M}_n)-E_+) \leq t+n^{-\epsilon} )+n^{-\frac{1}{6}+C\epsilon}.\nonumber
\end{align}
When $\gamma=1$, similar results also hold for $\lambda_p(\phi^{-1/2}\bm{M}_n)$.
\end{theorem}

This universality theorem shows that the asymptotic distributions of $\lambda_1(\bm{M}_n)$, $\lambda_p(\bm{M}_n)$, and consequently $\sigma_1(\bm{M}_n)=\|\bm{M}_n\|_{\text{op}}$, rely solely on the first two moments of $\bm{X}$. In contrast, \citet{qiu2023asymptotic} and \citet{chernozhuokov2022improved} show that other widely-used norms of $\bm{M}_n$, such as $\|\bm{M}_n\|_{\text{F}}$ and $\|\bm{M}_n\|_{\text{sup}}$, are influenced by the first four moments of $\bm{X}$.
This universal characteristic of the operator norm offers a straightforward yet effective framework for bootstrap. 
Before applying universality to the universal bootstrap, we give the following corollary of Theorem \ref{gauappbd}, which is of notable independent interest.

\begin{Co}\label{co1}
    Consider the case $\gamma>1$, $\bm{\Sigma}=\bm{I}_p$, $\bm{R}=\bm{0}$. Define $\mu=p$, $\sigma=p^{1/2}n^{-1/6}$. Under Assumption \ref{ass1} for $\bm{Z}=(\bm{Z}_1,\cdots,\bm{Z}_n)^T$, we have as $n\to \infty$,
    \begin{align}\label{tw1}
        \frac{\lambda_1(\bm{Z}^T\bm{Z})-\mu}{\sigma}\Rightarrow \text{TW}_1,
    \end{align}
    where ``$\Rightarrow$'' represents weak convergence, and $\text{TW}_1$ is the Tracy-Widom law of type $1$.
\end{Co}

One of the central problems in random matrix theory is establishing the asymptotic distribution of the largest eigenvalue of a sample covariance matrix.
In this context, (\ref{tw1}) was derived for Gaussian $\bm{Z}$ in the setting where $p/n \to \infty$ \citep{karoui2003largest}. However, for general distributions of $\bm{Z}$ where $\gamma > 1$, no further results have been established. Our universality Theorem \ref{gauappbd} addresses this significant gap in random matrix theory.

\subsection{Consistency of universal bootstrap}\label{covtest}
We shall now provide the consistency results for our universal bootstrap using universality results established in Section \ref{univer}.
We remark that the theorems derived thus far apply specifically to the extreme eigenvalues $\lambda_1(\bm{M}_n)$ and $\lambda_p(\bm{M}_n)$. As we will show, extending these results to extreme singular values $\sigma_1(\bm{M}_n)$ requires only weaker assumptions.

\begin{manualtheorem}{\ref*{ass4}$'$}\label{ass5} 
When $\gamma=1,$ we require that there exists $\tau>0$ such that
    
    (\rmnum{1}) if $|E_+|>|E_{-}|$, we assume $| 1+\phi^{-1/2}\tilde{m}(E_{+}) \lambda_i(\bm{\Sigma}(E_{+}))| \geq \tau$ for $i=1,\cdots,p$.

    (\rmnum{2}) if $|E_+|<|E_{-}|$, we assume $| 1+\phi^{-1/2}\tilde{m}(E_{-}) \lambda_i(\bm{\Sigma}(E_{-}))| \geq \tau$ for $i=1,\cdots,p$.

    (\rmnum{3}) if $|E_+|=|E_{-}|$, we assume $| 1+\phi^{-1/2}\tilde{m}(E_{\pm}) \lambda_i(\bm{\Sigma}(E_{\pm}))| \geq \tau$ for $i=1,\cdots,p$.
\end{manualtheorem}

Assumption \ref{ass5} is weaker than Assumption \ref{ass4}. While this relaxation appears minor,
consider the case $\gamma=1$ and $\bm{\Sigma}=-\bm{R}=\bm{I}_p$. We can show that $\bm{M}_n'$ satisfies Assumption \ref{ass5} for all $n$ and $p$, but fails to satisfy Assumption \ref{ass4} if $p>n$.
This demonstrates the relevance of introducing Assumption \ref{ass5} when analyzing extreme singular values.
Recalling that $T$ is the largest singular value of $\bm{\Sigma}-\bm{\Sigma}_0$, we are led to define the concept of an admissible pair.

\begin{D}[Admissible pair]
    We call two $p$ by $p$ non-negative matrices $(\bm{\Sigma},\bm{\Sigma}_0)$ an admissible pair, if $(\bm{\Sigma},\bm{R})=(\bm{\Sigma},-\bm{\Sigma}_0)$ satisfy Assumptions \ref{ass2}, \ref{ass3}, \ref{ass5}.
\end{D}

Equipped with these results, we present the universal consistency theorem for $T$.
To establish the theoretical validity of this universal bootstrap procedure, we provide a bound on the uniform Gaussian approximation error
\begin{align}\label{tvdif2}
    \rho_n(\bm{\Sigma}_0)=\sup_{t\geq 0 }\bigg|\mathbb{P}\bigg(\bm{\hat{\Sigma}}\in \mathbf{B}_{\text{op}}(\bm{\Sigma}_0,t) \bigg)-\mathbb{P}\bigg(\bm{\hat{\Sigma}}^{\text{ub}}\in \mathbf{B}_{\text{op}}(\bm{\Sigma}_0,t) \bigg)  \bigg|,
\end{align}
where $ \mathbf{B}_{\text{op}}(\bm{\Sigma}_0,t) = \big\{ \bm{N} \ :  \|\bm{N}-\bm{\Sigma}_0\|_{\text{op}} \leq t \big\}$ is the operator norm ball centered at $\bm{\Sigma}_0$ with radius $t$. The probability in (\ref{tvdif2}) is defined with respect to $\bm{X}$ and $\bm{Y}$, whose rows have covariance matrix $\bm{\Sigma}$, which may differ from $\bm{\Sigma}_0$.

\begin{theorem}[Universal bootstrap consistency]\label{thmubc}
     Under Assumption \ref{ass1}, for any admissible pair $(\bm{\Sigma},\bm{\Sigma}_0)$, we have the uniform Gaussian approximation bound for some constant $\delta>0$,
     \begin{align}\label{gab}
    \rho_n(\bm{\Sigma}_0) \lesssim n^{-\delta}\to 0 \quad \text{as} \quad n \to \infty.
\end{align}
Moreover, we have with probability approaching $1$,
\begin{align}\label{ubc}
    \sup_{0\leq\alpha\leq 1 }\bigg|\mathbb{P}\bigg(T \geq \hat{q}^{\text{ub},B}_{\bm{\Sigma},\bm{\Sigma}_0}(\alpha)\bigg|\bm{Y}^1,\cdots,\bm{Y}^B  \bigg)-\alpha \bigg| \lesssim \rho_n(\bm{\Sigma}_0) + B^{-\frac{1}{2}}.
\end{align}
\end{theorem}

Theorem \ref{thmubc} establishes the uniform Gaussian approximation bound and Type \uppercase\expandafter{\romannumeral1} error of the universal bootstrap.
Expression (\ref{ubc}) also ensures that we can uniformly control the test size with an error of at most $n^{-\delta}+B^{-1/2}$, which vanishes as $n\to \infty$ and $B\to\infty$. This result guarantees the uniform consistency of the universal bootstrap procedure. 

\begin{remark}
    To contextualize our findings, we briefly compare (\ref{gab}) with existing high-dimensional bootstrap results. Generally, these results are presented as
\begin{equation*}
    \rho(\mathcal{A})=\sup_{A\in \mathcal{A}}\big|\mathbb{P}\big(T\in A\big)-\mathbb{P}\big(T^*\in A \ \big| \ \bm{X}\big)  \big| \to 0,
\end{equation*}
where $T$ is a statistic, $T^*$ is its Gaussian counterpart, and $\mathcal{A}$ represents a specified family of sets. For example, \citet{chernozhukov2013gaussian, chernozhukov2017central, chernozhukov2023nearly} considered $T$ as mean estimators and $\mathcal{A}$ as all rectangular in $\mathbb{R}^p$. A more related choice of $\mathcal{A}$ is in \citet{zhai2018high,xu2019pearson,fang2024large}, who also considered mean estimators but take $\mathcal{A}$ to be the sets of Euclidean balls and convex sets in $\mathbb{R}^p$. Their results demonstrated that for sets of Euclidean balls $\mathcal{A}$, $\rho(\mathcal{A})$ converge to $0$ if and only if $p/n\to 0$, meaning the Gaussian approximation holds when $p=o(n)$. For comparison, we observe that the operator norm ball $\mathbf{B}_{\text{op}}(\bm{v},t)$ for vector $\bm{v}$ in $\mathbb{R}^p$ coincides with the Euclidean balls, and our results show that the universality approximation holds when $p/n$ converges to a nonzero constant or even diverges to infinity. For the covariance statistics, \citet{han2018gaussian} took $T$ as the sample covariance matrix and $\mathcal{A}$ as all sets of $s$-sparse operator norm balls (defined in their work). Especially, with $\mathcal{A}$ as all operator norm balls, i.e. $s=p$,  they required $p=o(n^{1/9})$, limiting $p$ to be considerably smaller than $n$.
Similarly, \citet{lopes2022improved,lopes2023bootstrapping} considered $T$ as sample covariance with $\mathcal{A}$ as operator norm balls, but imposed a decay rate $i^{-\beta}$ for the $i$-th largest eigenvalue $\lambda_i(\bm{\Sigma})$ of $\bm{\Sigma}$ with $\beta >1$, implying a low intrinsic test dimension. Likewise, \citet{giessing2023gaussian} required the effective rank $r(\bm{\Sigma})=\tr(\bm{\Sigma})/|\bm{\Sigma}|_{\op}$ to satisfy $r(\bm{\Sigma})=o(n^{1/6})$.
In contrast, we impose no such assumptions, allowing each eigenvalue of $\bm{\Sigma}$ to be of comparable scale.
These comparisons underscore the advantages of our proposal, even in regimes lacking consistency.
\end{remark}

\subsection{Consistency of generalized universal bootstrap}\label{geub}
In this section, we show that the universal bootstrap consistency extends to the general spectral statistic class (\ref{extrf}). An important application is the combined statistic $T^{\text{Com}}$ in (\ref{comb}).
While Theorem \ref{thmubc} provides thresholds for $T$ and $T^{\text{Roy}}$ individually, it cannot be directly applied to $T^{\text{Com}}$ due to the complex dependence between $T$ and $T^{\text{Roy}}$. To address this, we develop the generalized universality theorem.
We first introduce
\begin{align}\label{gopball}
    \mathbf{B}_{k,\text{op}}(\bm{\Sigma}_1,\bm{\Sigma}_2,\bm{t}_k) = \bigg\{ \bm{N} \ : \ \sigma_i\left(\bm{\Sigma}_2^{-\frac{1}{2}}(\bm{N}-\bm{\Sigma}_1) \bm{\Sigma}_2^{-\frac{1}{2}}\right) \leq t_i,\ i=1,\cdots,k \bigg\},
\end{align}
where $\bm{\Sigma}_2$ is positive-definite, and $\bm{t}_k=(t_1,\cdots,t_k)^T$ with $t_i\geq 0$ for each $i=1,\cdots,k$. Since the operator norm $\|\bm{N}-\bm{\Sigma}_1\|_{\text{op}}$ equals the largest singular value $\sigma_1(\bm{N}-\bm{\Sigma}_1)$, we obtain $\mathbf{B}_{\text{op}}(\bm{\Sigma},t)=\mathbf{B}_{1,\text{op}}(\bm{\Sigma},\bm{I}_p,t)$.
Next, we define the following quantity for symmetric matrices $\bm{\Sigma}_{1,m}$ and positive-definite matrices $\bm{\Sigma}_{2,m}$, $m=1,\cdots,M$,
\begin{align}
    \rho_n( \left\{ \bm{\Sigma}_{1,m},\bm{\Sigma}_{2,m}\right\}_{m=1}^M  )=\sup_{\bm{t}_{k_1},\cdots,\bm{t}_{k_M}\geq \bm{0} }\bigg|\mathbb{P}&\bigg(\bm{\hat{\Sigma}}\in \bigcap_{m=1}^M\mathbf{B}_{k_m,\text{op}}(\bm{\Sigma}_{1,m},\bm{\Sigma}_{2,m},\bm{t}_{k_m}) \bigg)\nonumber\\
    &-\mathbb{P}\bigg(\bm{\hat{\Sigma}}^{\text{ub}}\in \bigcap_{m=1}^M\mathbf{B}_{k_m,\text{op}}(\bm{\Sigma}_{1,m},\bm{\Sigma}_{2,m},\bm{t}_{k_m})  \bigg)  \bigg|,\label{gtvdif}
\end{align}
where we use $\bm{t}_k\geq \bm{0}$ to represent $t_i\geq 0$ for each $i=1,\cdots,k$.

\begin{theorem}[Generalized universal bootstrap consistency]\label{thmgubc}
     Under Assumption \ref{ass1}, given $M$ admissible pairs $(\bm{\Sigma}_{2,m}^{-1/2}\bm{\Sigma}\bm{\Sigma}_{2,m}^{-1/2},\bm{\Sigma}_{2,m}^{-1/2}\bm{\Sigma}_{1,m}\bm{\Sigma}_{2,m}^{-1/2})$ and fixed integer numbers $k_m$ for $m=1,\cdots,M$,
     we have the uniform Gaussian approximation bound for some constant $\delta>0$,
     \begin{align}\label{ggab}
    \rho_n( \left\{ \bm{\Sigma}_{1,m},\bm{\Sigma}_{2,m}\right\}_{m=1}^M  ) \lesssim n^{-\delta}\to 0 \quad \text{as} \quad n \to \infty.
\end{align}
Moreover, we have
\begin{align}\label{gubc}
    \sup_{0\leq\alpha\leq 1 }\bigg|\mathbb{P}\bigg(T^{\text{ExS}}\geq \hat{q}^{\text{ExS,ub,B}}_{\bm{\Sigma},\left\{\bm{\Sigma}_{1,m},\bm{\Sigma}_{2,m}\right\}_{m=1}^{M}}(\alpha) \bigg)-\alpha \bigg| \lesssim \rho_n(\bm{\Sigma}_0) + B^{-\frac{1}{2}}.
\end{align}
\end{theorem}

Theorem \ref{thmgubc} extends Theorem \ref{thmubc} by demonstrating that universal bootstrap consistency applies to general statistics based on extreme singular values with arbitrary dependencies in the form of (\ref{extrf}).  
This result reinforces our universal bootstrap principle: to approximate the distribution of statistics involving extreme eigenvalues, it suffices to substitute all random variables with their Gaussian counterparts. 

\section{Application to covariance inference}\label{apps}
In this section, we discuss applications of the universal bootstrap in high-dimensional covariance inference. We first apply our universal bootstrapped statistics to the covariance test $H_0: \bm{\Sigma}=\bm{\Sigma}_0$ with operator norm.
As noted in Section \ref{intro}, several popular statistics utilize the operator norm, including $T^{\text{Roy}}$ in (\ref{roy}). Our Theorem \ref{thmubc}, \ref{thmgubc} ensures the size control of universal bootstrap for these statistics. In this subsection, we further perform a power analysis across a family of statistics, with $T$ and $T^{\text{Roy}}$ as specific cases. 

We conduct a power analysis under the generalized spike model setting as \citet{jiang2021generalized}, which allows variability in the bulk eigenvalues and block dependence between the spiked and bulk parts.
We assume the following generalized spike structure
\begin{align}\label{gspk1}
    \bm{\Sigma}=\bm{U}\left[
	\begin{array}{cc}
		\bm{D} &\bm{0} \\
		\bm{0}	& \bm{V}_2
	\end{array}
	\right]\bm{U}^T,\ \bm{\Sigma}_0=\bm{U}\left[
	\begin{array}{cc}
		\bm{V}_1 &\bm{0} \\
		\bm{0}	& \bm{V}_2
	\end{array}
	\right]\bm{U}^T,\ \bm{\Sigma}_1=\bm{U}\left[
	\begin{array}{cc}
		\bm{R}_1 &\bm{0} \\
		\bm{0}	& \bm{R}_2
	\end{array}
	\right]\bm{U}^T,
\end{align}
where $\bm{D},\bm{V}_1,\bm{R}_1\in \mathbb{R}^{k\times k}$, $\bm{V}_2,\bm{R}_2\in\mathbb{R}^{(p-k)\times (p-k)}$ are diagonal matrices, $k$ is a fixed number, and $\bm{U}$ is orthogonal matrix. 
We consider the family of statistics $T(\bm{\Sigma}_1)=\left\| \bm{\Sigma}_1^{-\frac{1}{2}} \left(\hat{\bm{\Sigma}} -\bm{\Sigma}_0\right)  \bm{\Sigma}_1^{-\frac{1}{2}}\right\|_{\op}$,
and reject $H_0$ if $T(\bm{\Sigma}_1)\geq \hat{q}^{\text{ub},B}_{\bm{\Sigma}_1^{-1/2}\bm{\Sigma}_0\bm{\Sigma}_1^{-1/2},\bm{\Sigma}_1^{-1/2}\bm{\Sigma}_0\bm{\Sigma}_1^{-1/2}}(\alpha)$.
Notably, $T=T(\bm{I}_p)$, $T^{\text{Roy}}=T(\bm{\Sigma}_0)$ are both included in this family.
Define $\bm{D}'=\bm{D}\bm{R}^{-1}_1$, $\bm{V}_1'=\bm{V}_1\bm{R}^{-1}_1$, $\bm{V}_2'=\bm{V}_2\bm{R}^{-1}_2$.
We assume that $\|\bm{D}'\|_{\op}$, $\|\bm{V}'_1\|_{\op}$, $\|\bm{V}'_2\|_{\op}$ are bounded.
Define $E_+(\bm{\Sigma}_1)$, $E_{-}(\bm{\Sigma}_1)$ and $\tilde{m}(z;\bm{\Sigma}_1)$ as in Section \ref{univ} with $(\bm{\Sigma},\bm{R})=(\bm{V}_2',-\bm{V}_2')$ in $\bm{M}_n$. 
For simplicity, we only consider $\gamma=1$.

\begin{theorem}\label{power}
    Under Assumptions \ref{ass1}, define $\bm{D}'=\text{diag}\left((d_i')_{i=1}^{k}\right)$, $\bm{V}_1'=\text{diag}\left((v_i')_{i=1}^{k}\right)$ and 
    \begin{align}\label{kappap}
        \kappa= -\frac{\phi^{\frac{1}{2}}}{\tilde{m}(\phi^{-\frac{1}{2}}E_+(\bm{\Sigma}_1);\bm{\Sigma}_1)}-\frac{\phi^{\frac{1}{2}}v_i'}{E_+(\bm{\Sigma}_1)\tilde{m}(\phi^{-\frac{1}{2}}E_+(\bm{\Sigma}_1);\bm{\Sigma}_1)} .
    \end{align}

    (\rmnum{1}) If $d_i' >\kappa$ for some $1\leq i \leq k$, then the power goes to $1$, i.e.
    \begin{align*}
        \mathbb{P}\left(T(\bm{\Sigma}_1)\geq \hat{q}^{\text{ub},B}_{\bm{\Sigma}_1^{-1/2}\bm{\Sigma}_0\bm{\Sigma}_1^{-1/2},\bm{\Sigma}_1^{-1/2}\bm{\Sigma}_0\bm{\Sigma}_1^{-1/2}}(\alpha) \right) \to 1,\ \text{as}\ n\to \infty, \ B\to\infty.
    \end{align*}

    (\rmnum{2}) If $d_i' <\kappa$ for all $1\leq i \leq k$ and $T(\bm{\Sigma}_1)=\lambda_1\left( \bm{\Sigma}_1^{-1/2} \left(\hat{\bm{\Sigma}} -\bm{\Sigma}_0\right)  \bm{\Sigma}_1^{-1/2}\right)$, then
    \begin{align*}
        \mathbb{P}\left(T(\bm{\Sigma}_1)\geq \hat{q}^{\text{ub},B}_{\bm{\Sigma}_1^{-1/2}\bm{\Sigma}_0\bm{\Sigma}_1^{-1/2},\bm{\Sigma}_1^{-1/2}\bm{\Sigma}_0\bm{\Sigma}_1^{-1/2}}(\alpha) \right) \to \alpha,\ \text{as}\ n\to \infty, \ B\to\infty,
    \end{align*}
    i.e. there is no power under this setting.
\end{theorem}

Theorem \ref{power} provides the power analysis of the test based on the statistic $T(\bm{\Sigma}_1)$. Specifically, it identifies the sharp phase transition point $\kappa$ in (\ref{kappap}),
which closely resembles the well-known Baik-Ben Arous-P\'{e}ch\'{e} (BBP) phase transition for the largest eigenvalues of the sample covariance matrix, see \citet{baik2005phase,baik2006eigenvalues,paul2007asymptotics}. When the spike eigenvalues exceed this phase transition point, they lie outside the bulk eigenvalue support, yielding full test power. Conversely, if the spike eigenvalues $d'_i$ fall below this point, distinguishing them from the null case using $T(\bm{\Sigma}_1)$ becomes infeasible.

With Theorem \ref{power}, we can now deduce the power of $T$ and $T^{\text{Roy}}$ as special cases.

\begin{Co}\label{compp}
    Under assumptions of Theorem \ref{power}, define $\bm{D}=\text{diag}\left((d_i)_{i=1}^{k}\right)$, $\bm{V}_1=\text{diag}\left((v_i)_{i=1}^{k}\right)$.

    (\rmnum{1}) If $ d_i -v_i > -\frac{\phi^{1/2}}{\tilde{m}(\phi^{-1/2}E_+)}-v_i\left(1+\frac{\phi^{1/2}}{E_+\tilde{m}(\phi^{-1/2}E_+)}\right)$ for some $1\leq i \leq k$,
    then the power of $T$ goes to $1$. Here $\phi^{-1/2}E_+\tilde{m}(\phi^{-1/2}E_+)<-1$.

    (\rmnum{2}) If $d_i -v_i >  \phi^{1/2} v_i$ for some $1\leq i \leq k$, then the power of $T^{\text{Roy}}$ goes to $1$.
\end{Co}

We compare the phase transition points (\rmnum{1}) and (\rmnum{2}) for $T$ and $T^{\text{Roy}}$, respectively. These points reveal distinct behaviors: as $\phi^{-1/2}E_+\tilde{m}(\phi^{-1/2}E_+)<-1$, the phase transition point of (\rmnum{1}) decreases with $v_i$, whereas the phase transition point of (\rmnum{2}) increases with $v_i$. This indicates that when spikes affect the larger eigenvalues, i.e. for large $v_i$, the statistics $T$ achieves high power more effectively than the normalized $T^{\text{Roy}}$ with a lower phase transition point. Conversely, the normalized statistic $T^{\text{Roy}}$ performs better for smaller $v_i$.

As another important application of the universal bootstrap, we show how to construct sharp uniform simultaneous confidence intervals for $\langle \bm{A}, \bm{\Sigma} \rangle$ for all $\bm{A} \in \mathbb{R}^{p \times p}$, where the inner product of two matrices $\bm{A}$ and $\bm{B}$ is defined as $\langle \bm{A}, \bm{B} \rangle = \tr(\bm{A}^T \bm{B})$. Notably, setting $\bm{A} = \bm{c}_1^{}\bm{c}_1^{T}$ and $\bm{A} = \bm{c}_1^{}\bm{c}_2^{T}$ gives the variance $\text{Var}(\bm{c}_1^{T}\bm{X}_i)$ and covariance $\text{Cov}(\bm{c}_1^{T}\bm{X}_i, \bm{c}_2^{T}\bm{X}_i)$ of all linear combinations $\bm{c}_1^{T}\bm{X}_i$ and $\bm{c}_2^{T}\bm{X}_i$ of $\bm{X}_i$, which has significant applications in optimal portfolio analysis \citep{bai2010spectral, hult2012risk}.
The main challenge comes from that the full matrix $\bm{\Sigma}$ and consequently $\hat{q}^{\text{ub},B}_{\bm{\Sigma}, \bm{\Sigma}}(\alpha)$ cannot be consistently estimated in our setting, so the universal bootstrap cannot be directly applied. To address this, we note that $\hat{q}^{\text{ub},B}_{\bm{\Sigma},\bm{\Sigma}}(\alpha)=\hat{q}^{\text{ub},B}_{\text{Sp}(\bm{\Sigma}),\text{Sp}(\bm{\Sigma})}(\alpha)$ where the spectral matrix is defined as $\text{Sp}(\bm{\Sigma})=\text{diag}(\lambda_1(\bm{\Sigma}),\cdots,\lambda_p(\bm{\Sigma}))$, and certain estimators $\widehat{\text{Sp}}(\bm{\Sigma})$ for the spectral matrix $\text{Sp}(\bm{\Sigma})$ are consistent under a suitable normalized distance. For instance, given the distance $d^2_{\text{Sp}}(\widehat{\text{Sp}}(\bm{\Sigma}),\text{Sp}(\bm{\Sigma}))=\tr (\widehat{\text{Sp}}(\bm{\Sigma})-\text{Sp}(\bm{\Sigma}))^2/p$, \citet{ledoit2015spectrum,kong2017spectrum} proposed estimators $\widehat{\text{Sp}}(\bm{\Sigma})$ satisfying $d_{\text{Sp}}(\widehat{\text{Sp}}(\bm{\Sigma}),\text{Sp}(\bm{\Sigma}))\to 0$ as $n\to \infty$ when $\gamma=1$. 

\begin{theorem}\label{usci}
    Under Assumption \ref{ass1}, for any admissible pair $(\bm{\Sigma},\bm{\Sigma})$ and estimators $\widehat{\text{Sp}}(\bm{\Sigma})$ of $\text{Sp}(\bm{\Sigma})$, we have the following sharp uniform simultaneous confidence intervals for $\langle \bm{A},\bm{\Sigma} \rangle$
    \begin{align}\label{usci1}
        \sup_{0\leq\alpha\leq 1 }\bigg|\mathbb{P}\bigg(\forall \bm{A}\in\mathbb{R}^{p\times p},\ \langle \bm{A},\bm{\Sigma} \rangle \in \bigg[\langle \bm{A},\bm{\hat{\Sigma}} \rangle-\hat{q}^{\text{ub},B}_{\widehat{\text{Sp}}(\bm{\Sigma}),\widehat{\text{Sp}}(\bm{\Sigma})}(\alpha)\ &\|\bm{A}\|_{S_{1}},\langle \bm{A}, \bm{\hat{\Sigma}} \rangle+\\
         \hat{q}^{\text{ub},B}_{\widehat{\text{Sp}}(\bm{\Sigma}),\widehat{\text{Sp}}(\bm{\Sigma})}(\alpha)\ \|\bm{A}\|_{S_{1}}\bigg] \bigg|\bm{Y}^1,\cdots,\bm{Y}^B\bigg)  \nonumber  
        -(1-\alpha)        \bigg| 
         \lesssim \rho_n(\bm{\Sigma}_0)+&d_{\text{Sp}}(\widehat{\text{Sp}}(\bm{\Sigma}),\text{Sp}(\bm{\Sigma})) + B^{-\frac{1}{2}}.
    \end{align}
    Here the Schatten $1-$norm of matrix $\bm{A}$ is defined as $\|\bm{A}\|_{S_{1}}=\sum_{i=1}^{p}\sigma_{i}(\bm{A})$.
    As a special case, we also have simultaneous confidence intervals for $\bm{c}_1^T \bm{\Sigma} \bm{c}_2$
    \begin{align}\label{usci2}
        \sup_{0\leq\alpha\leq 1 }\bigg|\mathbb{P}\bigg(\forall \bm{c}_1,\bm{c}_2\in \mathbb{R}^p,\ \bm{c}_1^T \bm{\Sigma} \bm{c}_2 \in \bigg[\bm{c}_1^T \bm{\hat{\Sigma}} \bm{c}_2-\hat{q}^{\text{ub},B}_{\widehat{\text{Sp}}(\bm{\Sigma}),\widehat{\text{Sp}}(\bm{\Sigma})}&(\alpha)\ \|\bm{c}_1\|\cdot\|\bm{c}_2\|,\ \bm{c}_1^T \bm{\hat{\Sigma}} \bm{c}_2+\\
        \hat{q}^{\text{ub},B}_{\widehat{\text{Sp}}(\bm{\Sigma}),\widehat{\text{Sp}}(\bm{\Sigma})}(\alpha)\ \|\bm{c}_1\|\cdot\|\bm{c}_2\| \bigg]\bigg|\bm{Y}^1,\cdots,\bm{Y}^B\bigg)  \nonumber  
        -(1-\alpha) \bigg| 
        &\lesssim \rho_n(\bm{\Sigma}_0) +d_{\text{Sp}}(\widehat{\text{Sp}}(\bm{\Sigma}),\text{Sp}(\bm{\Sigma}))+ B^{-\frac{1}{2}}. 
    \end{align}
\end{theorem}

Theorem \ref{usci} establishes sharp uniform simultaneous confidence intervals for both $\langle \bm{A},\bm{\Sigma} \rangle$ and $\bm{c}_1^T \bm{\Sigma} \bm{c}_2$.
Specifically, intervals (\ref{usci1}) and (\ref{usci2}) are sharp, as their probability converges uniformly to $1 - \alpha$, rather than no less than $1-\alpha$. Constructing confidence intervals for all linear combinations of $\bm{\Sigma}$ is considerably more challenging than for individual entries. For the mean vector case, achieving confidence intervals for each entry requires only $\ln p =o(n^a)$ for some $a>0$, see \citet{chernozhukov2013gaussian, chernozhukov2017central, chernozhukov2023nearly}. In contrast, confidence intervals for all linear combinations of the mean vector require $p=o(n)$, see \citet{zhai2018high,xu2019pearson,fang2024large}. Remarkably, we construct sharp simultaneous confidence intervals for all $\langle \bm{A},\bm{\Sigma} \rangle$ even as $p/n$ converges or diverges to infinity.

\section{Numerical results}\label{numericalres}
\subsection{Simulation}\label{simu}
In this subsection, we numerically analyze the performance of the universal bootstrap method in high-dimensional covariance tests. We evaluate the empirical size and power of the universal bootstrapped operator norm-based statistics $T$ in (\ref{statistic}) (Opn) and the combined statistics $T^{\text{Com}}$ in (\ref{comb}) (Com). These are compared with the operator norm-based statistic $T^{\text{Roy}}$ in (\ref{roy}) (Roy), the Frobenius norm-based linear spectral statistic $T^{\text{F},1}$ from \citet{qiu2023asymptotic} (Lfn), the debiased Frobenius norm-based U statistic $T^{\text{F},2}$ from \citet{cai2013optimal} (Ufn), and the supremum norm-based statistic $T^{\text{sup}}$ inspired by \citet{cai2011limiting,cai2013two} (Supn). The expressions for these statistics are provided in the Supplement. 

To evaluate the robustness of our method, we generate $n$ independent, non-identically distributed samples of $p$-dimensional random vectors with zero means and covariance matrix $\bm{\Sigma}$. We set sample sizes $n=100$ and $n=300$ and vary the dimension $p$ over $100$, $200$, $500$, $1000$, and $2000$ to cover both proportional and ultra-high-dimensional cases. For each configuration, we calculate the average number of rejections based on $2000$ replications, with the universal bootstrap procedure using $1000$ resamplings. The significance level is fixed at $0.05$.  We consider three different structures for the null covariance matrix $\bm{\Sigma}_0=(\sigma_{0,ij})_{p\times p}$:

(a). \textit{Exponential decay.} The elements $\sigma_{0,ij}=0.6^{|i-j|}$. 

(b). \textit{Block diagonal.} Let $K=\lfloor p/10\rfloor$ be the number of blocks, with each block being $1$ on the diagonal and $0.55$ on the off-diagonal. 

(c). \textit{Signed sub-exponential decay.} The elements $\sigma_{0,ij}=(-1)^{i+j}0.4^{|i-j|^{0.5}}$.

These structures have been previously analyzed in \citet{cai2013two,zhu2017testing,yu2024fisher} for models (a) and (b), and in \citet{cai2013two,yu2024fisher} for model (c), highlighting the broad applicability of our approach. To assess the empirical size, we examine the following distributions: \textit{Gaussian distribution}, \textit{Uniform and t-distributions}, \textit{Gaussian and uniform distributions}, and \textit{Gaussian and t-distributions}.
Here the \textit{Uniform and t-distributions} sets $\{Z_{ij} ,\ i=1,\cdots,\lfloor n/2\rfloor; j =1,\cdots,p\}$ to be i.i.d. normalized uniform random variables on $[-1,1]$, and $\{Z_{ij} ,\ i=\lfloor n/2\rfloor+1,\cdots,n; j =1,\cdots,p\}$ to be i.i.d. normalized t random variables with degree of freedom $12$.
The other distributions are defined similarly.
All distributions of $Z_{ij}$ are standardized to zero mean and unit variance.
Common covariance tests often consider the Gaussian and uniform distributions \citep{chen2023sharp}, while the t-distribution with $12$ degrees of freedom is examined in \citet{cai2013two,zhu2017testing}.  Except for the Gaussian baseline, the other cases share identical covariance matrices but differ in distribution, creating challenging scenarios for testing covariance. These reflect the broad applicability of our universality result.

Table \ref{size1} reports the empirical sizes of the supremum, the Frobenius, and the operator norm tests at the significance level $0.05$,  sample size $n=300$ across various dimension $p$, covariance $\bm{\Sigma}_0$ and data distributions of Gaussian and the uniform and t-distributions. 
For results on $n=100$ and for Gaussian-uniform and Gaussian-$t$ distributions with $n=100,300$, see Table 5-10 in the Supplement. The three operator norm tests are performed using the proposed universal bootstrap procedure, while the supremum and Frobenius norm tests rely on their respective asymptotic distribution formulas. 
Both Table \ref{size1} and Table 5-10 in the Supplement show that the universal bootstrap effectively controls the size of all operator norm-based statistics at the nominal significance level across all tested scenarios.

\begin{landscape}
\begin{table}
\centering
\caption{The empirical sizes of different norm tests with a significance level $0.05$, the sample size $n=300$ for data with Gaussian distribution and the t and uniform distribution and the combinations of the dimension $p$, and the covariance matrix $\bm{\Sigma}_0$.}
\label{size1}
    \begin{tabular}{ccccc ccccc ccccc cc}
        \hline \multirow{2}{*}{} & \multirow{2}{*}{\begin{tabular}{c}
                $\bm{\Sigma}_0$ \\
                $p$
        \end{tabular}}  & \multicolumn{5}{c}{ \textbf{Exp. decay(0.6)} } & \multicolumn{5}{c}{ \textbf{Block diagonal} }& \multicolumn{5}{c}{ \textbf{Signed subExp. decay(0.4)} } \\
        \cmidrule(r){3-7}\cmidrule(lr){8-12}\cmidrule(lr){13-17}
        & & \textbf{100} & \textbf{200} & \textbf{500} & \textbf{1000} &\textbf{2000} & \textbf{100} & \textbf{200} & \textbf{500} & \textbf{1000}&\textbf{2000} & \textbf{100} & \textbf{200} & \textbf{500} & \textbf{1000} &\textbf{2000}\\
        \hline & & \multicolumn{15}{c}{ Gaussian, n=300 } \\
        \multirow{7}{*}{  }
        &Supn  & 0.073 & 0.066 & 0.071&0.070 & 0.076& 0.073 &0.066 &0.071 &0.070  & 0.076 &0.073 &0.066 &0.071 &0.070&0.076\\
	& Lfn & 0.056 & 0.048 & 0.045&0.045 &0.053 & 0.056 &0.048 &0.045 &0.045  & 0.053 & 0.056& 0.048&0.045 &0.045&0.053\\
	&Ufn  & 0.053 & 0.048 & 0.044&0.044 & 0.057& 0.053 &0.048 &0.044 & 0.044 & 0.057  &0.053 &0.048 &0.044 &0.044&0.057\\
	&Roy  & 0.048 & 0.049 &0.050 & 0.052&0.051 & 0.048 &0.049 &0.050 & 0.052 &0.051  & 0.048&0.049 & 0.050&0.052&0.051\\
    & Opn & 0.057 & 0.052 &0.047 &0.051 &0.052 &0.052  &0.050 &0.052 & 0.052 & 0.049  & 0.042&0.046 &0.054 &0.049&0.054\\
	&Com  &0.057  & 0.050 &0.046 &0.052 &0.052 & 0.053 & 0.050& 0.050& 0.051 & 0.048  &0.042 &0.048 &0.052 &0.049&0.055\\
            \cmidrule(r){2-17}
            & & \multicolumn{15}{c}{  uniform and t(df=12), n=300 } \\
            \multirow{7}{*}{  }
           &Supn  &0.093  & 0.089 &0.098 &0.109 &0.125 & 0.093 &0.089 &0.098 &0.109  &0.125  & 0.093&0.089 &0.098 &0.109&0.125\\
		& Lfn & 0.055 & 0.047 &0.050 &0.039 &0.053 & 0.055 &0.047 &0.050 &0.039  & 0.053 &0.055 &0.047 &0.050 &0.039&0.053\\
		& Ufn  &0.056  & 0.044 &0.049 &0.040 &0.056 &0.056  &0.044 & 0.049& 0.040 & 0.056 &0.056 &0.044 &0.049 &0.040&0.056\\
		&Roy  & 0.050 & 0.046 &0.038 &0.053 &0.054 & 0.050 &0.046 &0.038 &0.053  & 0.054 &0.050 &0.046 &0.038 &0.053&0.054\\
        & Opn & 0.040 & 0.051 &0.049 &0.051 &0.050 & 0.042 & 0.053& 0.051&0.047  & 0.051 &0.046 & 0.042&0.036 &0.042&0.061\\
		&Com  & 0.046 & 0.051 & 0.041&0.048 &0.047 & 0.043 &0.052 &0.050 &0.046  & 0.052 & 0.049&0.045 &0.038 &0.045&0.058\\
            \hline
        \end{tabular}
\end{table}

\end{landscape}

\noindent  
This approach performs well under both Gaussian and non-Gaussian distributions, for both i.i.d. and non-i.i.d. data, in proportional and ultra-high dimensional settings, and across various covariance structures. These findings confirm the consistency of the universal bootstrap procedure, providing empirical support for our theoretical results in Theorem \ref{thmubc} and Theorem \ref{thmgubc}. Additionally, the Frobenius norm-based test maintains appropriate size control, while the supremum norm-based test exhibits substantial size inflation at $n=100$ and moderate inflation even at a larger sample size of $n=300$.

To evaluate the empirical power of these statistics, we consider a setting where $\bm{\Sigma}=\bm{\Sigma}_0+\bm{\Delta}$, with $\bm{\Delta}$ having a low rank.  This setup corresponds to the power analysis presented in Theorem \ref{power}. As shown in Theorem \ref{power}, the statistic $T$ outperforms $T^{\text{Roy}}$ when the eigenspaces of $\bm{\Delta}$ align with the eigenvectors of $\bm{\Sigma}_0$ corresponding to larger eigenvalues. In contrast, $T^{\text{Roy}}$ performs better when $\bm{\Delta}$ aligns with the smaller eigenvectors of $\bm{\Sigma}_0$. For a fair comparison, we define $\bm{\Delta}=\sigma\bm{u}\bm{u}^T/2+\sigma\bm{v}\bm{v}^T/4$, where $\bm{u}$ is the fifth-largest eigenvector of $\bm{\Sigma}_0$, $\bm{v}$ is uniformly sampled on the unit sphere in $\mathbb{R}^p$, and $\sigma$ is the signal strength. We term this configuration the \textit{spike setting}.
For completeness, we also conduct a power analysis where $\bm{\Delta}$ has full rank by setting $\bm{\Delta} = \sigma \bm{I}_{p}$, where $\sigma$ represents the signal level. This configuration, termed the \textit{white noise setting}, represents the covariance structure after adding white noise to the null covariance. Together, the spike and white noise settings cover both low- and full-rank deviations from the null covariance. 
For illustration, we consider sample sizes $n=100$ and $300$, with dimension $p=1000$. Based on the universality results in Theorem \ref{thmubc} and empirical size performance, we conduct the power analysis using a Gaussian distribution for simplicity. Further empirical power comparisons under varying sample sizes and dimensions are provided in Table 1-4 of the Supplement table.

\begin{figure}[!htb]
  \centering
  \subfloat[Exp. decay, n=100]{
    \includegraphics[width=4.5cm]{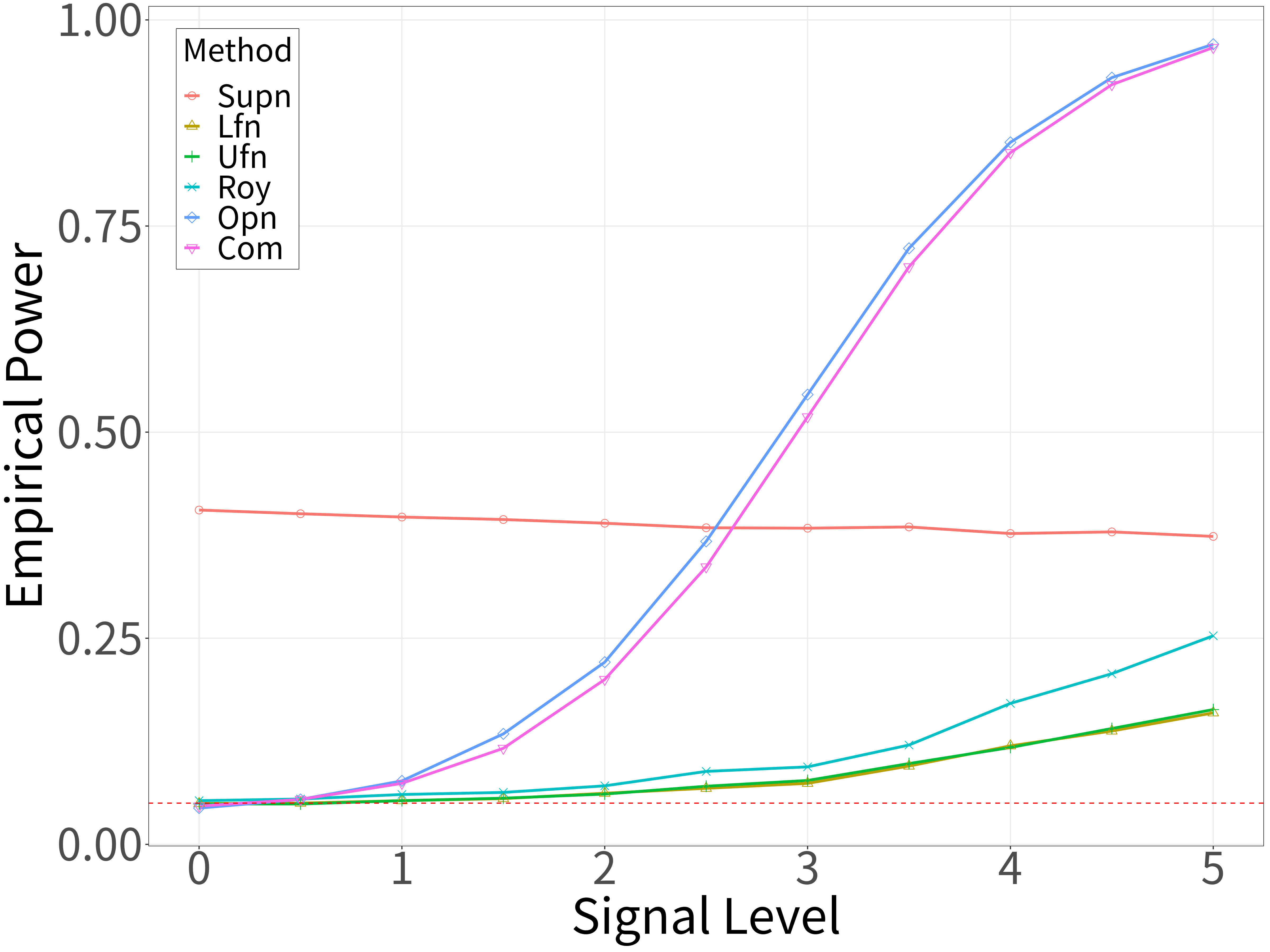}
  }
  \subfloat[Block diagonal, n=100]
  {
    \includegraphics[width=4.5cm]{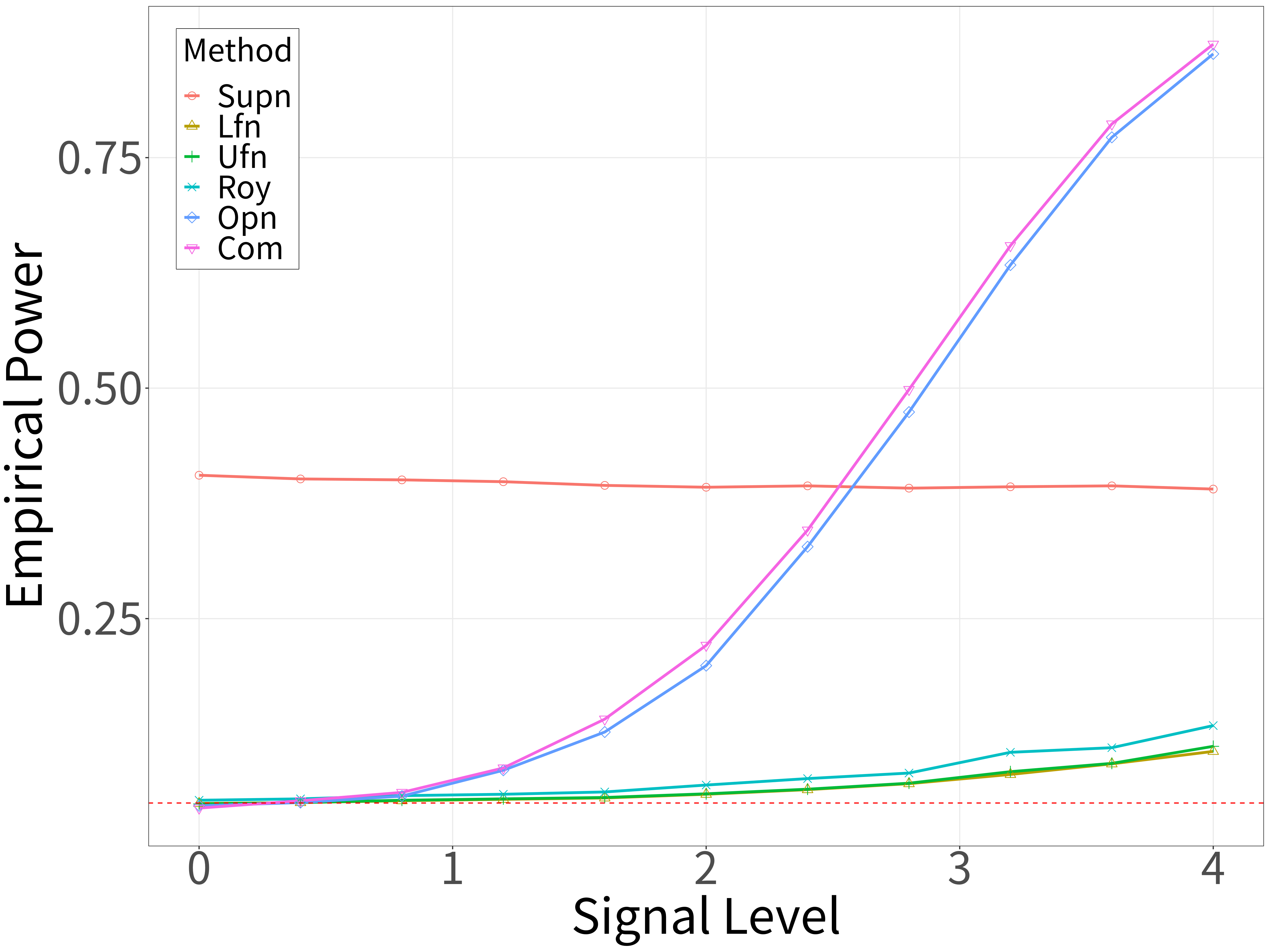}
  }
  \subfloat[Signed subExp, n=100]
  {
    \includegraphics[width=4.5cm]{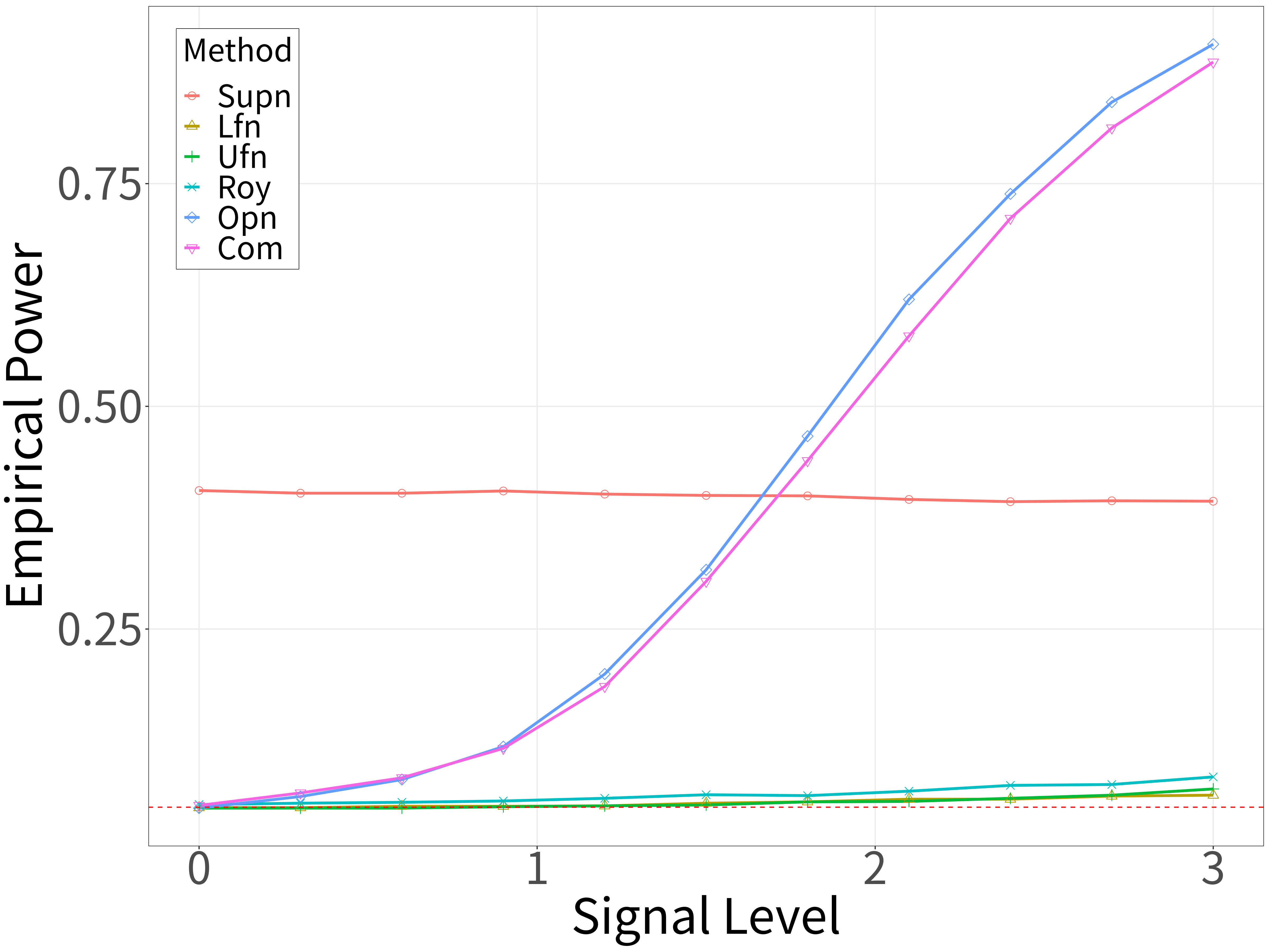}
  }
  
  \subfloat[Exp. decay, n=300]{
    \includegraphics[width=4.5cm]{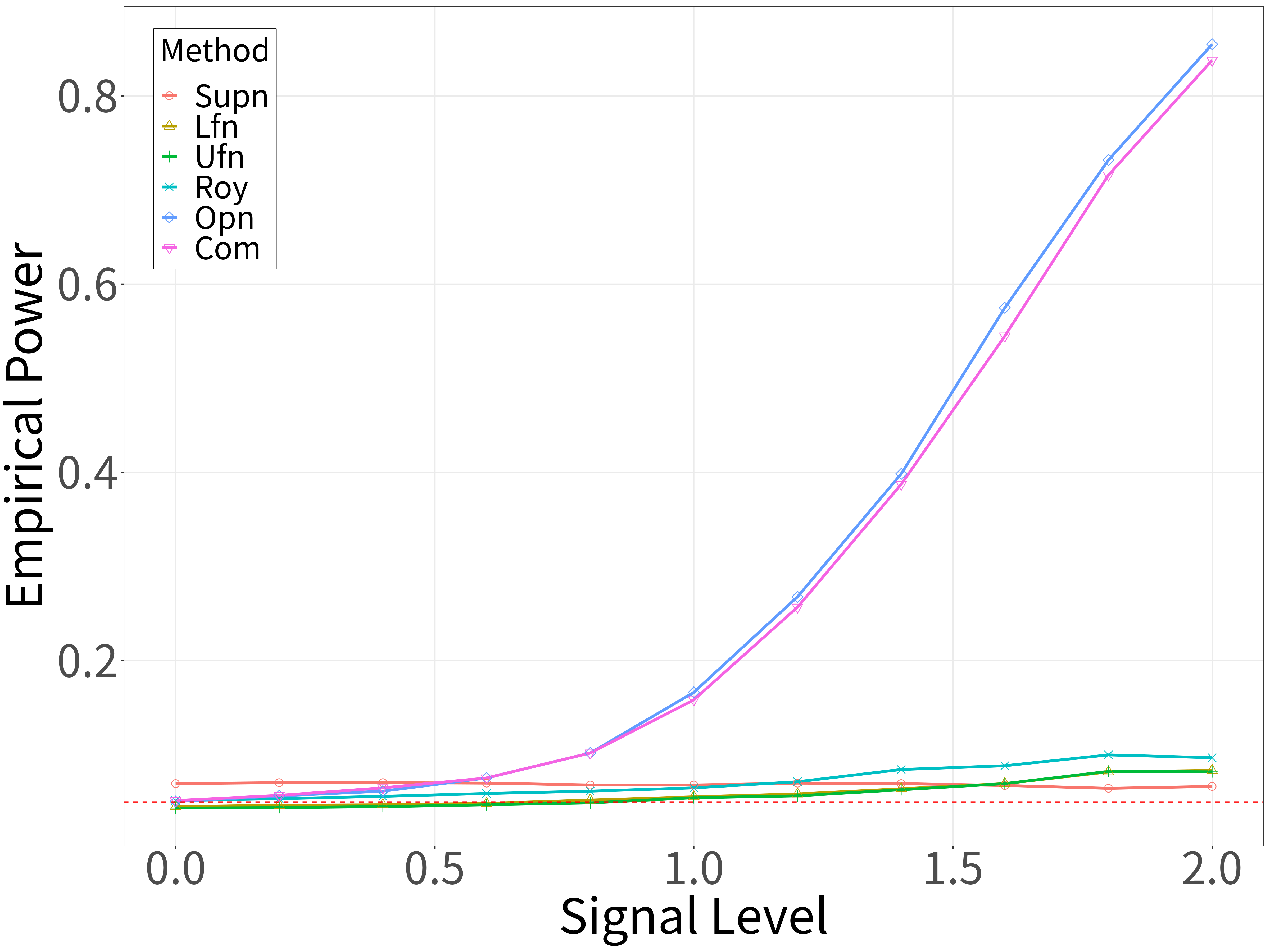}
  }
  \subfloat[Block diagonal, n=300]
  {
    \includegraphics[width=4.5cm]{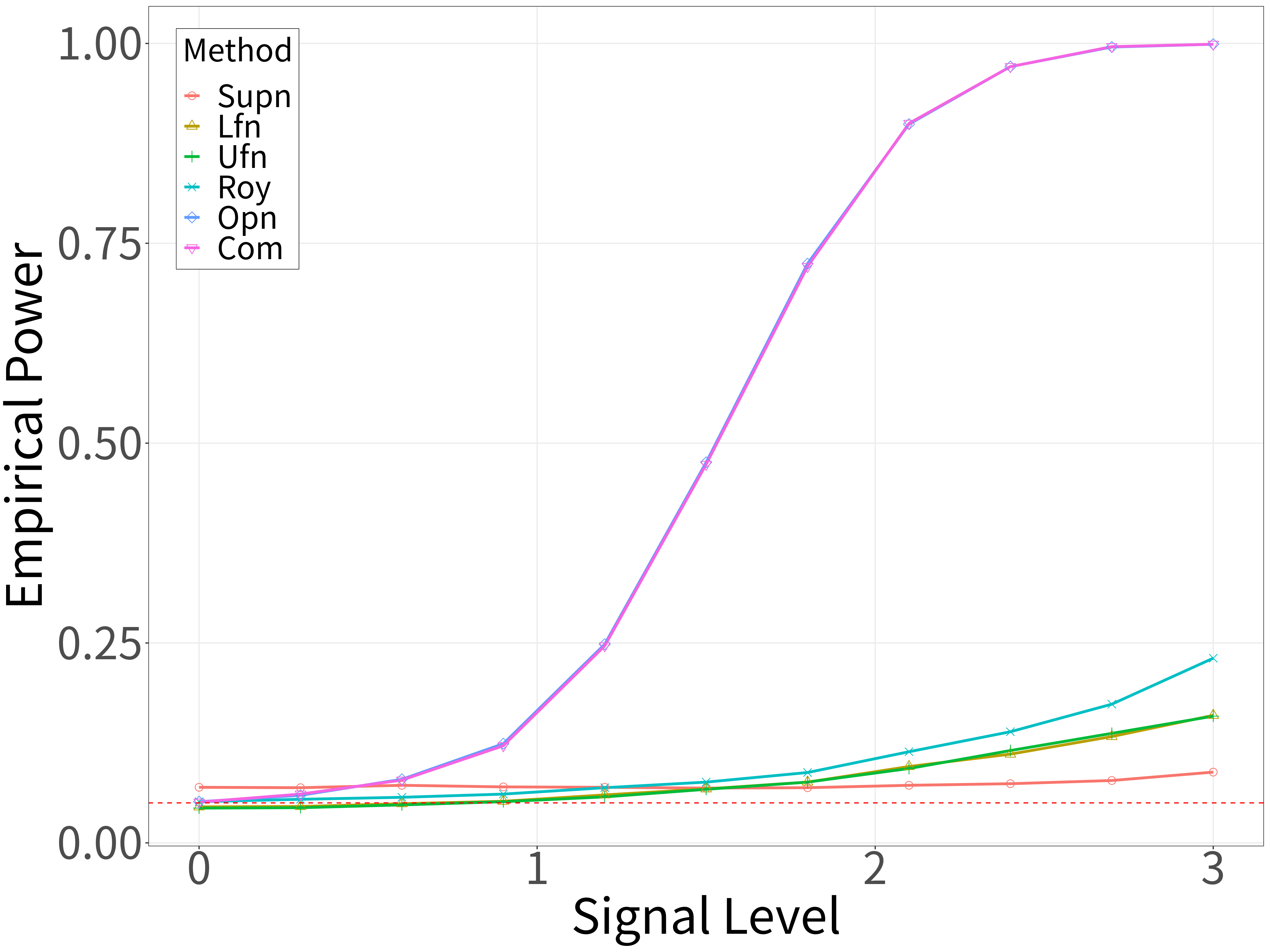}
  }
  \subfloat[Signed subExp, n=300]
  {
    \includegraphics[width=4.5cm]{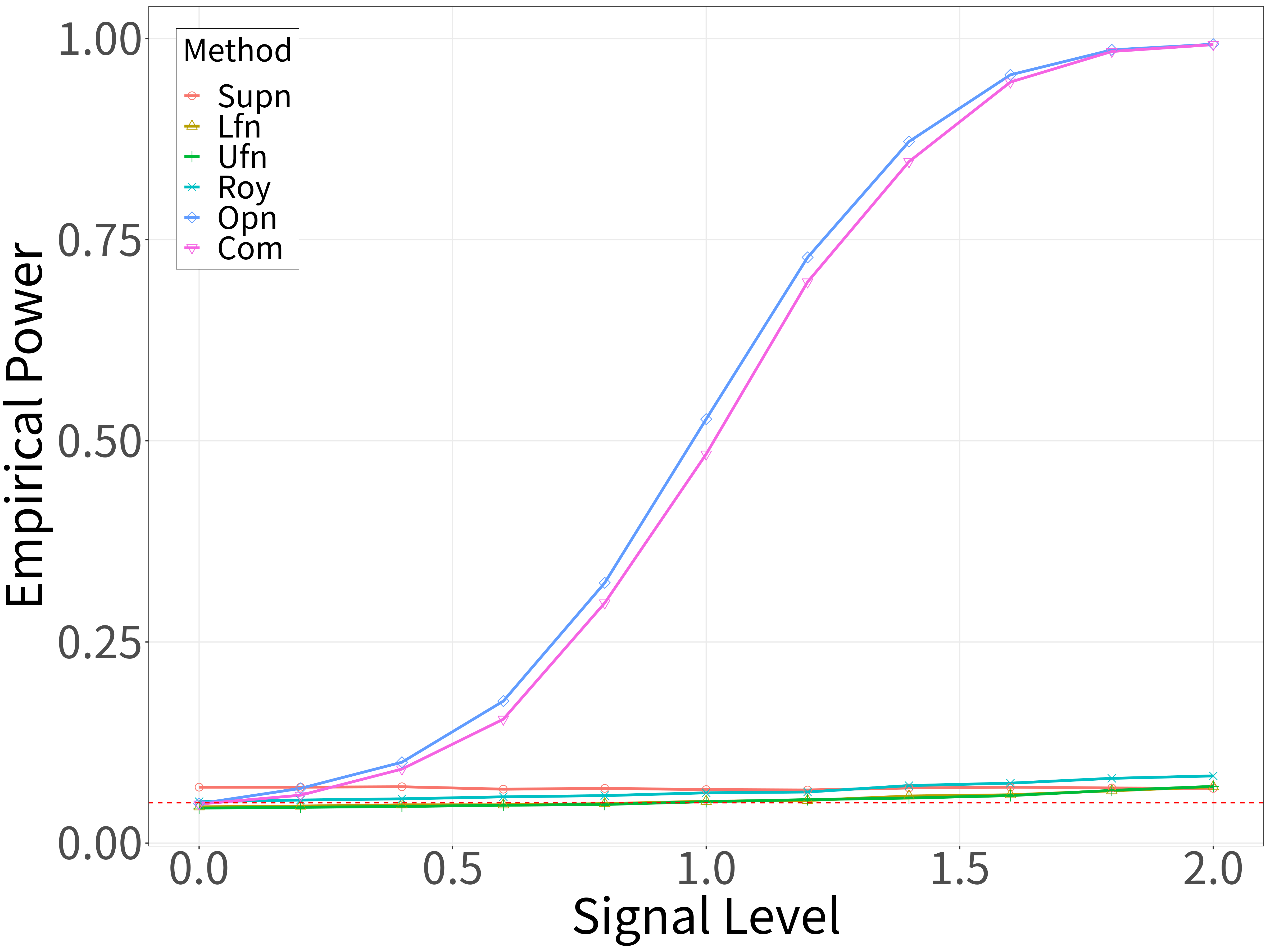}
  }
  \caption{Empirical powers of the supremum, the Frobenius, and the proposed operator norm tests with
respect to the signal level of the spike setting alternatives under three covariance structures, the sample size $n=100$, $300$, dimension $p=1000$, and the Gaussian data with $2000$ replications. }
  \label{fig1}
\end{figure}

\begin{figure}[!htb]
  \centering
  \subfloat[Exp. decay, n=100]{
    \includegraphics[width=4.5cm]{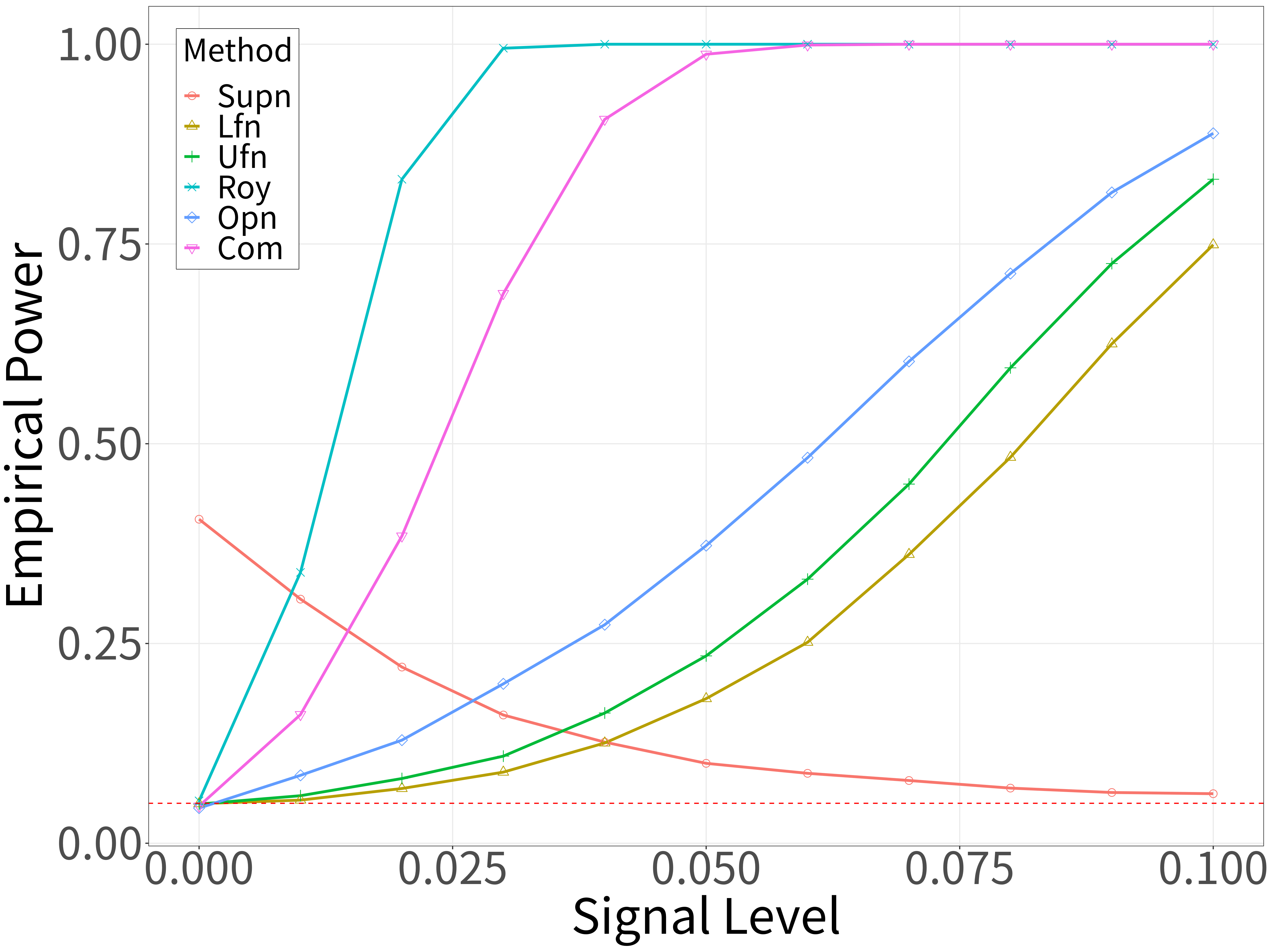}
  }
  \subfloat[Block diagonal, n=100]
  {
    \includegraphics[width=4.5cm]{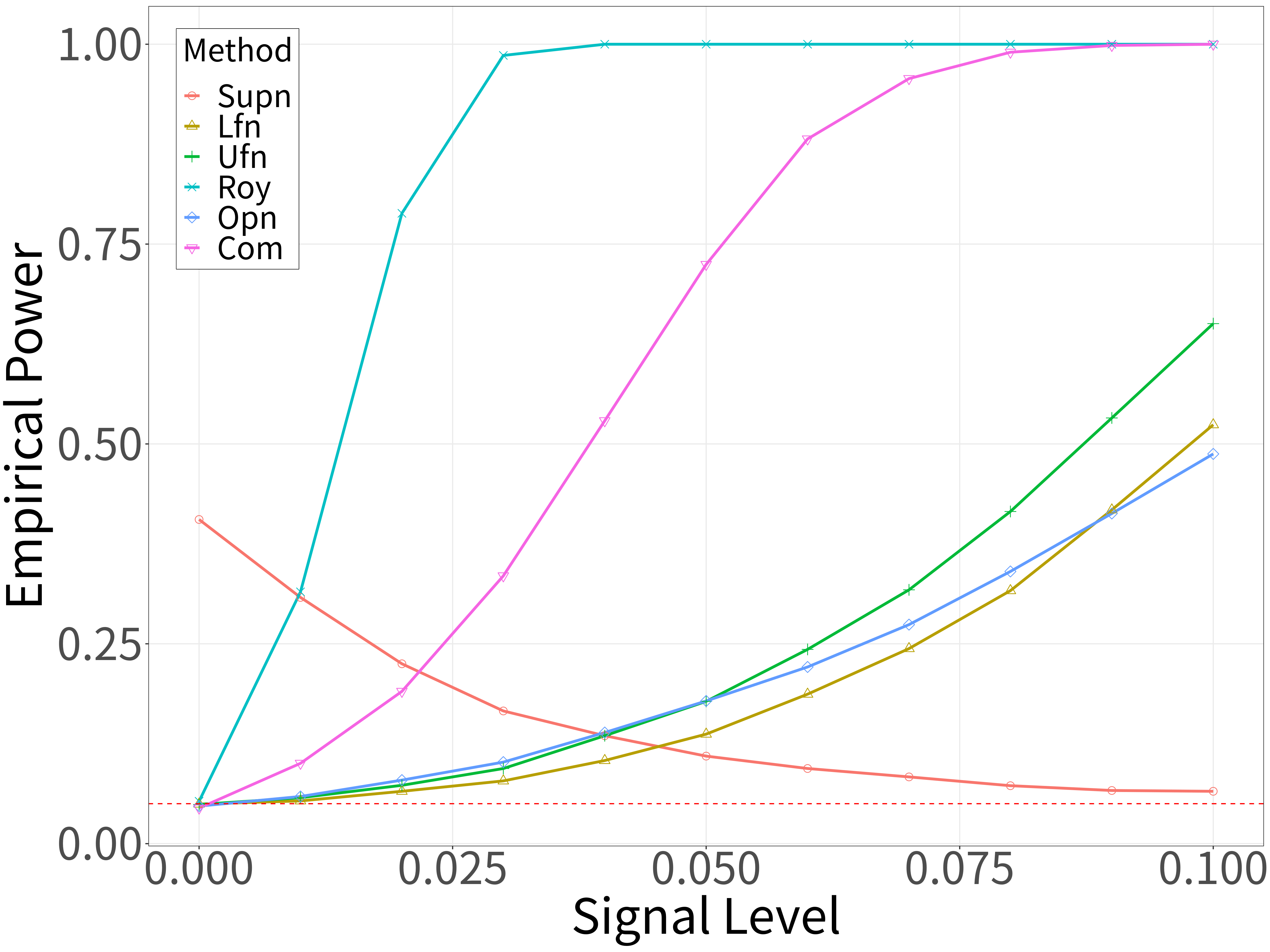}
  }
  \subfloat[Signed subExp, n=100]
  {
    \includegraphics[width=4.5cm]{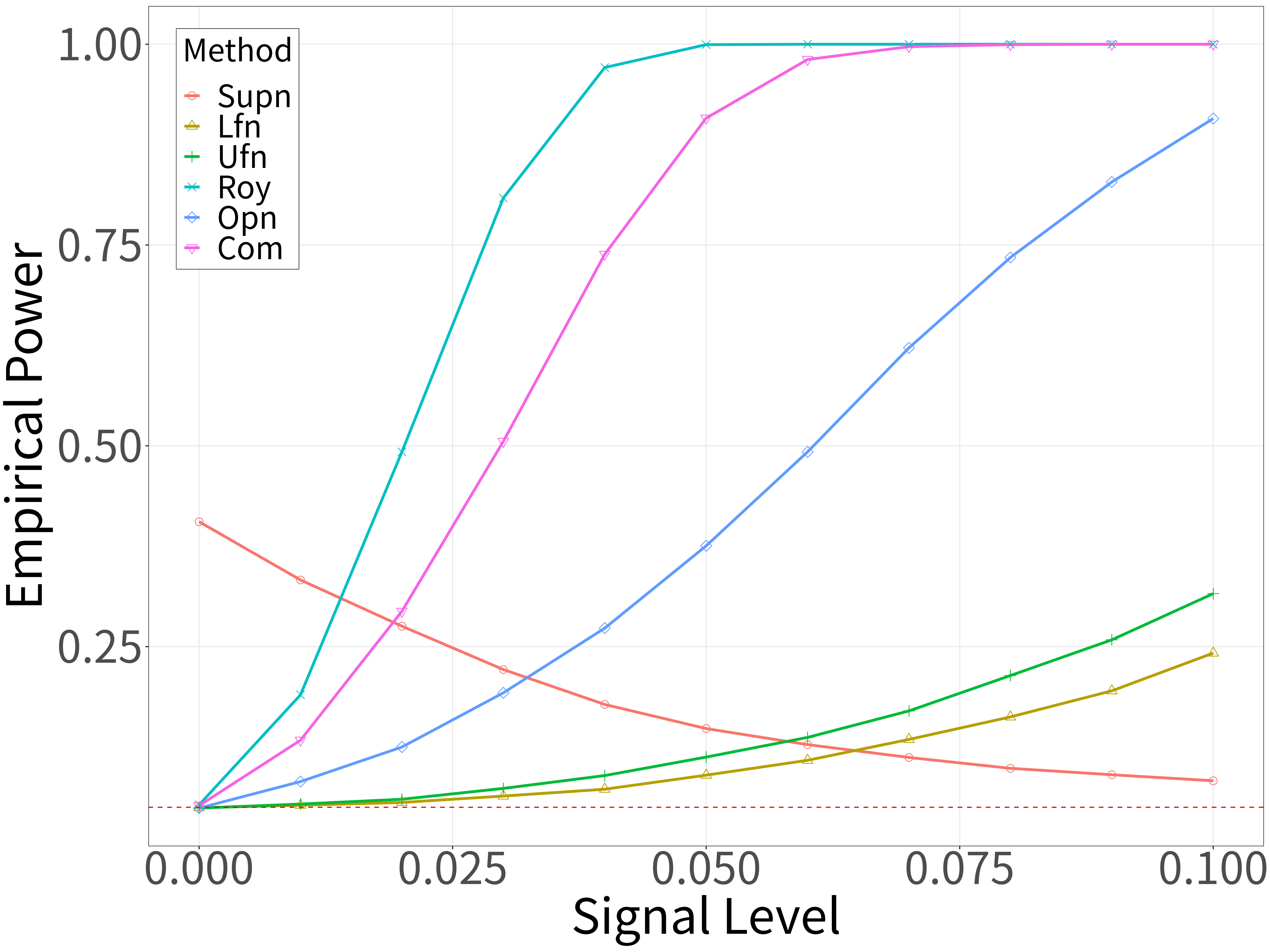}
  }
  
  \subfloat[Exp. decay, n=300]{
    \includegraphics[width=4.5cm]{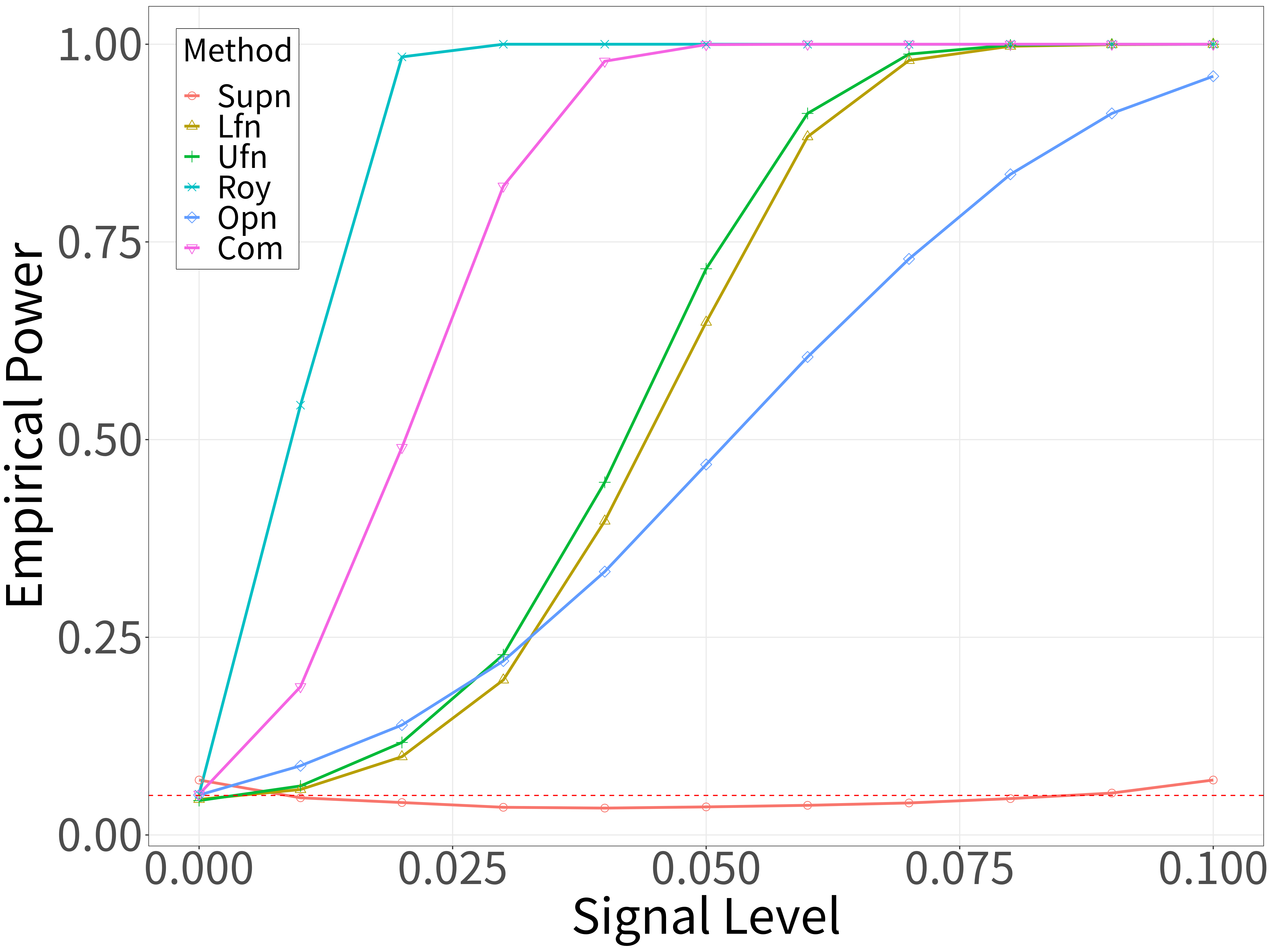}
  }
  \subfloat[Block diagonal, n=300]
  {
    \includegraphics[width=4.5cm]{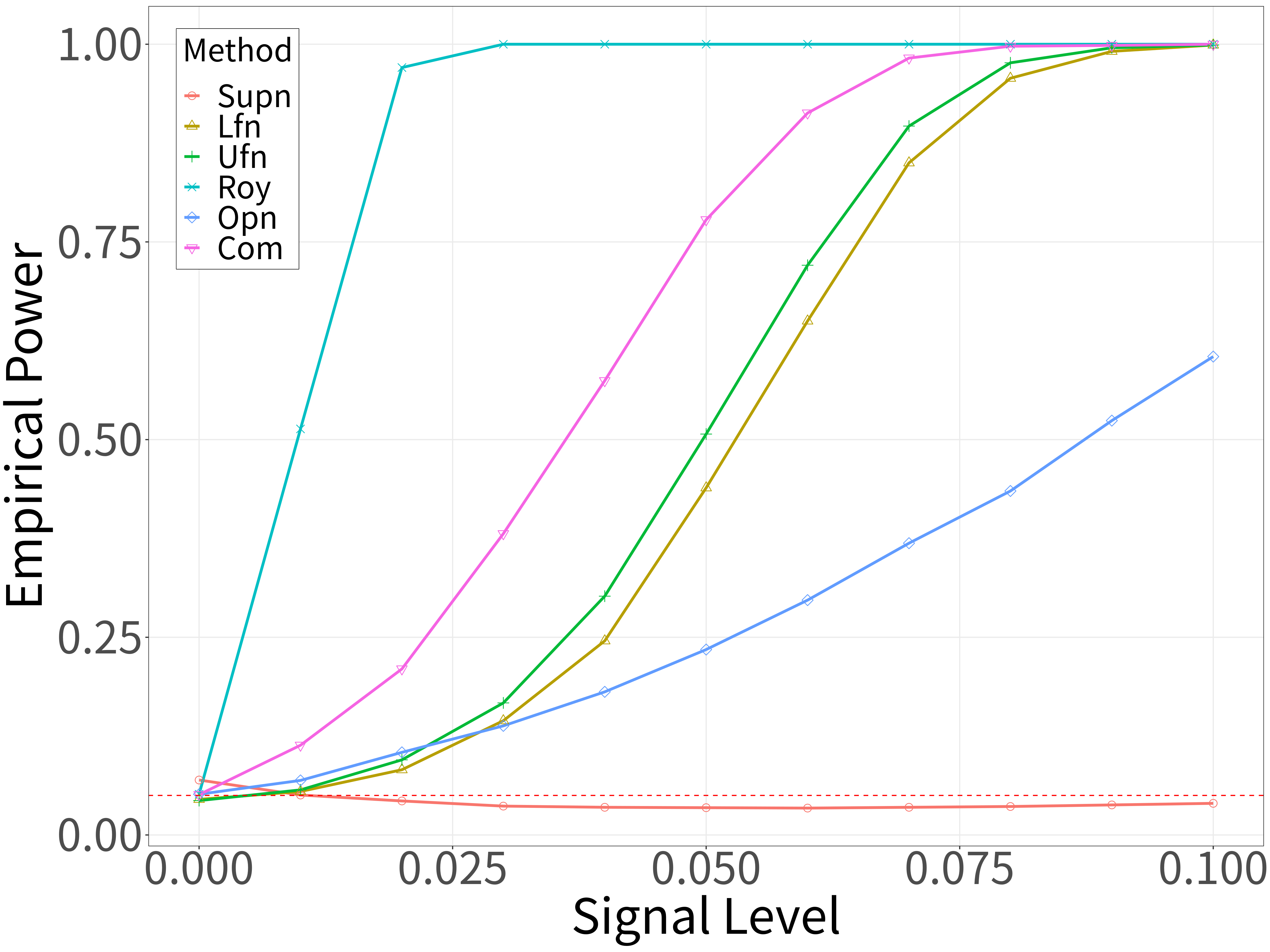}
  }
  \subfloat[Signed subExp, n=300]
  {
    \includegraphics[width=4.5cm]{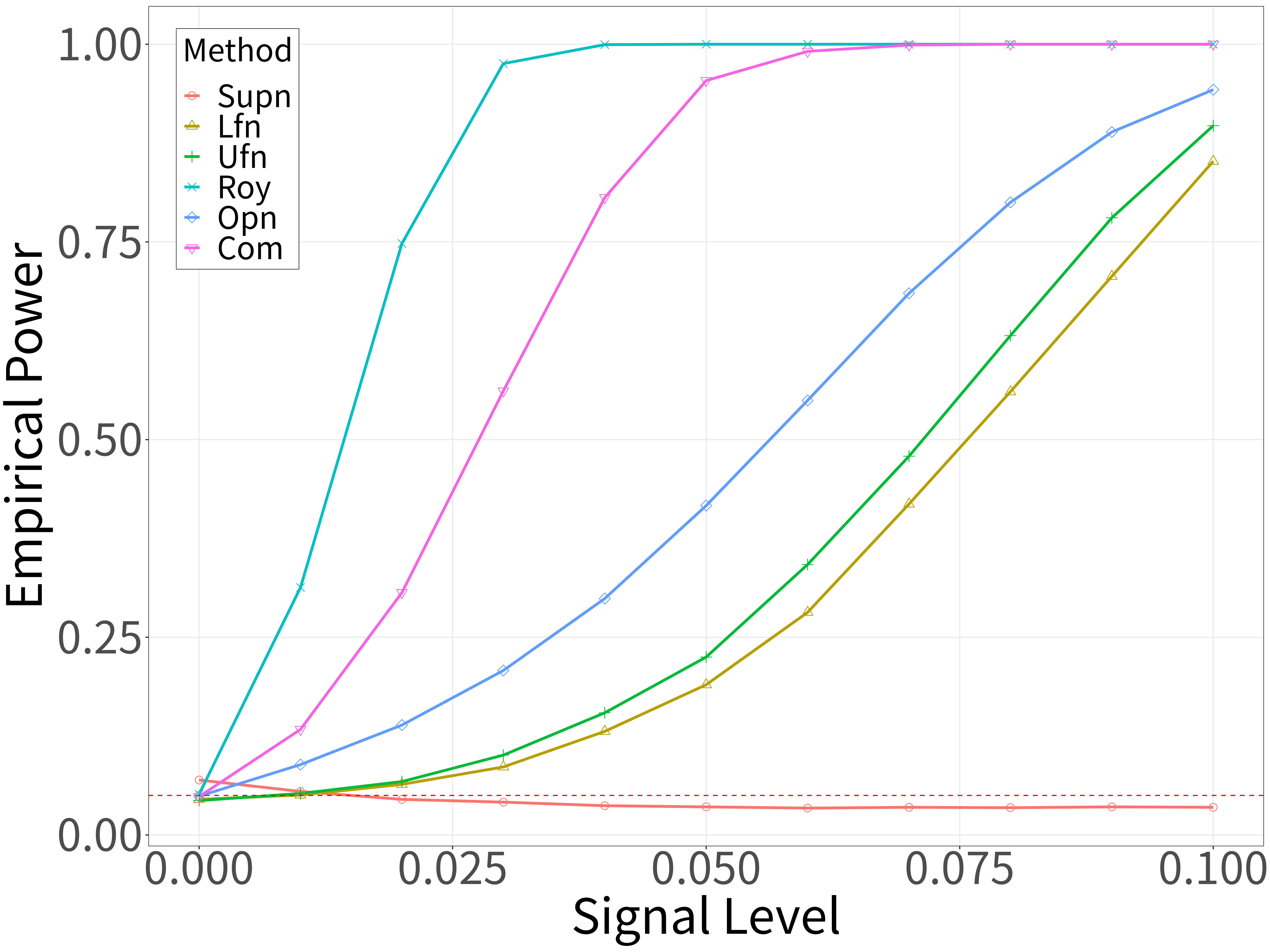}
  }
  \caption{Empirical powers of the supremum, the Frobenius, and the proposed operator norm tests with
respect to the signal level of the white noise setting alternatives under three covariance structures, the sample size $n=100$, $300$, dimension $p=1000$, and the Gaussian data with $2000$ replications.}
  \label{fig2}
\end{figure}

Figure \ref{fig1} illustrates the empirical power performance of six test statistics across different signal levels and null covariance structures under the spike setting, where the dashed line represents the significance level of $0.05$. Under this setting, the proposed operator norm-based statistics, $T$ and $T^{\text{Com}}$ exhibit similar power, while the two Frobenius norm-based statistics and the normalized operator norm statistic $T^{\text{Roy}}$ perform comparably. 
When $n=100$, the operator norm-based statistics $T$ and $T^{\text{Com}}$ maintain appropriate test size and achieve the highest power, approaching $1$ quickly across all three covariance configurations. The normalized operator norm statistic $T^{\text{Roy}}$ and the two Frobenius norm-based statistics $T^{\text{F},1}$ and $T^{\text{F},2}$ also control the test size but exhibit lower power, among which $T^{\text{Roy}}$ slightly outperforms.
The supremum norm-based statistic, however, exhibits significant size inflation and fails to detect signals, as its power curve does not increase with signal level. The pattern is similar for $n=300$, except that the power curves for the supremum norm-based statistic remain near the dashed line, indicating limited power performance. This provides empirical evidence for the advantage of operator norm-based statistics over the other norms in the spike setting, supporting the proposed statistic $T$ over $T^{\text{Roy}}$. 

We further make comparisons under the white noise setting across various signal levels in Figure \ref{fig2}. In this setting, the combined statistic $T^{\text{Com}}$ and the normalized statistic $T^{\text{Roy}}$ outperform other statistics. Because the white noise setting applies a uniform signal across eigenvector directions of $\bm{\Sigma}_0$ with both large and small eigenvalues, the normalized statistic $T^{\text{Roy}}$ achieves higher power than $T$, which is consistent with Theorem \ref{power}. Notably, the combined statistic $T^{\text{Com}}$ performs robustly in both the spike and white noise settings, demonstrating enhanced power without size inflation, as supported by the generalized universal bootstrap consistency Theorem \ref{thmgubc}. Similar to the spike setting, Frobenius norm-based statistics display lower power, while the supremum norm-based statistic fails to control size. In summary, when the covariance structure is known, we recommend using the statistic $T$. When the covariance structure is unknown, we prefer the more robust combined statistic, $T^{\text{Com}}$.

\subsection{Data application}\label{realdata}
We apply our test procedure to annual mean near-surface air temperature data from the global scale for the period $1960-2010$, using the HadCRUT4 dataset \citep{morice2012quantifying} and the Coupled Model Intercomparison Project Phase 5 \citep[CMIP5,][]{taylor2012overview}, as detailed in \citet{li2021uncertainty} (\url{https://figshare.com/articles/dataset/Uncertainty_in_Optimal_Fingerprinting_is_Underestimated/14981241/3?file=28924416}). The global dataset includes monthly anomalies of near-surface air temperature across $5^{\circ}\times 5^{\circ}$ grid boxes. To reduce dimensionality, these grid boxes are aggregated into larger $40^{\circ}\times 30^{\circ}$ boxes, resulting in $S=54$ spatial grid boxes. 

To mitigate distribution shifts due to long-term trends, we divide the $1960-2010$ period into five decades and analyze the data for each decade separately. The monthly data for each decade are averaged over five-year intervals, yielding a temporal dimension $T=2$ for each ten years. We first outline the commonly used modeling procedure in climatology. Following \citet{li2021uncertainty}, one considers a high-dimensional linear regression
\begin{align*}
    \bm{Y}_k=\sum_{i=1}^2\bm{X}_{k,i}  \beta_{k,i} +\bm{\epsilon}_k,\ k=1,\cdots,5,
\end{align*}
where $\bm{Y}_k \in \mathbb{R}^{L}$ is the vectorized observed mean temperature across spatial and temporal dimensions for the $k$-th decades after $1960$, with dimension $L=S\times T$. The vectors $\bm{X}_{k,1}$, $\bm{X}_{k,2},\bm{\epsilon}_k \in \mathbb{R}^{L}$ represent the unobserved expected climate response under the ANT and NAT forcings, and the unobserved noise, respectively. The coefficients $\beta_{k,1}$, $\beta_{k,2}$ are unknown scaling factors of interest for statistical inference.
To estimate $\beta_{k,1}$, $\beta_{k,2}$, in addition to the observed data $\bm{Y}_k$, we have $n_i$ noisy observation fingerprints $\tilde{\bm{X}}_{k,i,j}$ of the climate system
\begin{align*}
    \tilde{\bm{X}}_{k,i,j}=\bm{X}_{k,i}+\tilde{\bm{\epsilon}}_{k,i,j},\ k=1,\cdots,5;\ j=1,\cdots,n_i;\ i=1,2.
\end{align*}
In our dataset, $n_1 = 35$, $n_2 = 46$.
For privacy, the data $\tilde{\bm{X}}_{k,i,j}$ are preprocessed by adding white noise, with the according covariance of the noise added to the hypothesized matrix $\bm{\Sigma}_{\bm{\epsilon},k}$ defined below.
The climate models also provide $N$ simulations $\hat{\bm{\epsilon}}_{k,i}$, $k=1,\cdots,5;\ i=1,\cdots,N$ with $N=223$ provided in this dataset, which are assumed to follow the same distribution as the unobserved noise $\bm{\epsilon}_k$. To efficiently estimate $\beta_1,\beta_2$, a typical assumption is that the natural variability simulated by the climate models matches the observed variability, i.e., the covariance of $\tilde{\bm{\epsilon}}_{k,i,j}$, $k=1, \dots, 5; \ j=1, \dots, n_i; \ i=1,2$, equals the covariance matrix $\bm{\Sigma}_{\bm{\epsilon},k}$ of $\bm{\epsilon}_k$. Since $\bm{\epsilon}_k$ is unobserved, a common choice for $\bm{\Sigma}_{\bm{\epsilon},k}$ is the sample covariance of the simulations $\bm{\Sigma}_{\bm{\epsilon},k}=\sum_{i=1}^{N}\hat{\bm{\epsilon}}_{k,i}^{}\hat{\bm{\epsilon}}_{k,i}^T/N$. Therefore, it is important to test the hypothesis
\begin{align}\label{exinapp}
    H_{0,k} : \text{Cov}(\tilde{\bm{\epsilon}}_{k,i,j})=\bm{\Sigma}_{\bm{\epsilon},k},\ k=1,\cdots,5;\ j=1,\cdots,n_i;\ i=1,2.
\end{align}
As noted by \citet{olonscheck2017consistently}, this equivalence (\ref{exinapp}) is crucial for optimal fingerprint estimation in climate models. However, few studies have validated (\ref{exinapp}) with statistical evidence, which calls for the need to test this hypothesis.

\begin{table}
\centering
\caption{Estimated p-values for observation group and control group using the supremum norm-based statistic (Supn), the Frobenius norm-based statistic (Lfn, Ufn), and the operator norm-based statistics (Roy, Opn, Com), with sample size $n=81$, dimension $p=108$.}
\label{realta}
\begin{tabular}{ccccccccccccc}
\hline p-value & \multicolumn{6}{c}{ observation group } & \multicolumn{6}{c}{ control group } \\
\cmidrule(r){2-7}\cmidrule(lr){8-13}
year & Supn  &Lfn &Ufn &Roy & Opn & Com& Supn  &Lfn &Ufn &Roy & Opn & Com  \\
\hline 
1960-1970 & .783 &.006&.001 & .000 &.000  & .000 &.002&.566 & .521 &.805 &.646&.734\\
1970-1980 & .406 &.146&.115& .000 & .000 & .000 &.004&.230&.243  & .401&.452&.427 \\
1980-1990 & .249 &.001&.000& .000 &  .000& .000 &.000&.419& .451 &.638&.713&.677 \\
1990-2000 & .905 &.634&.410& .000 & .000 &.000 &.223&.072&.042  &.445 &.383&.414 \\
2000-2010 & .984 &.002&.001& .000 & .000 & .000 &.719&.368&.376  & .244 &.288&.268\\
\hline
\end{tabular}
\end{table}

We construct the data matrix $\tilde{\bm{X}}_k \in \mathbb{R}^{n \times p}$ by combining quantities $\tilde{\bm{X}}_{k,1,j} - \bar{\tilde{\bm{X}}}_{k,1}$ for $j = 1, \dots, n_1$ and $\tilde{\bm{X}}_{k,2,j} - \bar{\tilde{\bm{X}}}_{k,2}$ for $j = 1, \dots, n_2$, where $n = n_1 + n_2 = 81$ and $p = L = 108$.
Here $\bar{\tilde{\bm{X}}}_{k,i}=\sum_{j=1}^{n_i}\tilde{\bm{X}}_{k,i,j}/n_i$. We apply several statistics to the data matrix $\tilde{\bm{X}}_k$ with the hypothesized matrix $\bm{\Sigma}_{\bm{\epsilon},k}$: the proposed operator norm-based statistics $T$ (Opn) and  $T^{\text{Com}}$ (Com), the operator norm-based statistics $T^{\text{Roy}}$ (Roy), the Frobenius norm-based statistics $T^{\text{F},1}$ (Lfn) and $T^{\text{F},2}$ (Ufn), and the supremum norm-based statistic $T^{\text{sup}}$ (Supn). For comparison, we generate i.i.d. Gaussian samples $\tilde{\bm{X}}_{k,j}' \sim \mathcal{N}(\bm{0}, \bm{\Sigma}_{\bm{\epsilon},k})$ for $j = 1, \dots, n;\ k = 1, \dots, 5$ as a control group. The observed data $\tilde{\bm{X}}_k$ forms the observation group.
The p-value results for these statistics are summarized in Table \ref{realta}. The supremum norm-based statistic $T^{\text{sup}}$ fails to reject the observation group for all periods but rejects the control group in the periods $1960-1970$, $1970-1980$, and $1980-1990$. The Frobenius norm-based statistics $T^{\text{F},1}$ and $T^{\text{F},2}$ fail to reject the observation group in the periods $1970-1980$ and $1990-2000$, and $T^{\text{F},2}$ rejects the control group during $1990-2000$. Only the operator norm-based statistics $T$, $T^{\text{Com}}$, and $T^{\text{Roy}}$ reject the null hypothesis in the observation group while not rejecting the control group for all years. This suggests that the commonly assumed hypothesis (\ref{exinapp}) should be rejected and replaced by a more suitable assumption to estimate $\beta_{k,1}$ and $\beta_{k,2}$.

\bigskip
\begin{center}
{\large\bf SUPPLEMENTARY MATERIAL}
\end{center}

\begin{description}

\item[Universal\_Bootstrap\_supp:] The additional theoretical results related to Section \ref{univ} and \ref{apps}, more simulation examples and implementation details, all proofs and auxiliary lemmas are deferred to the supplementary material. (.pdf type)

\item[Code and Data] Code to implement and reproduce the simulations and real data applications and corresponding output and raw datasets. 

\end{description}

\bibliographystyle{agsm}

\bibliography{Bibliography-MM-MC}
\end{document}